\documentclass[12pt]{amsart}
\usepackage{brownian}

\togglefalse{EA} 

\begin{document}

\begin{abstract}
A Brownian motion tree (BMT) model is a Gaussian model whose associated set of covariance matrices is linearly constrained according to common ancestry in a phylogenetic tree. 
We study the complexity of inferring the maximum likelihood (ML) estimator for a BMT model by computing its ML-degree.
Our main result is that the ML-degree of the BMT model on a star tree with $n + 1$ leaves is $2^{n+1}-2n-3$, which was previously conjectured by Am\'endola and Zwiernik.  We also prove that the ML-degree of a BMT model is independent of the choice of the root.
The proofs rely on the toric geometry of concentration matrices in a BMT model. Toward this end, we
produce a combinatorial formula for the determinant of the concentration matrix of a BMT model, which generalizes the Cayley-Pr\"ufer theorem to complete graphs with weights given by a tree. 
\end{abstract}
\maketitle
\section{Introduction}
\label{sec:Introduction}

A Brownian motion tree (BMT) model is a family of multivariate Gaussian distributions that describe the evolution of a continuous trait in a set of $n$ species along a phylogenetic tree. 
{In this setting, a rooted phylogenetic tree on $n+1$ leaves is a tree with no degree 2 vertices and with leaves $0,\dots,n$ where all edges of the tree are directed away from leaf $0$ \cite{boege2021reciprocal,felsenstein1973maximum, SUZ20}.}
Non-root leaves correspond to the extant species of interest, and the other nodes in the tree are the common ancestors of these species.
{In many contexts, the root of a phylogenetic tree is considered to be a degree two vertex which all edges are directed away from. However, in the study of Brownian motion tree models, it is instead convenient to consider the root to be a \emph{leaf} labeled $0$, which can be viewed as a common ancestor of the taxa or an outgroup species \cite{felsenstein1973maximum}.}

Introduced by Felsenstein \cite{felsenstein1973maximum} in 1973, BMT models enjoy many applications in phylogenetics. They have been applied to test for selective pressure \cite{cooper2010body,freckleton2006detecting}, often by serving as null model for evolution under genetic drift \cite{schraiber2013inferring}. BMT models are commonly used to represent continuous molecular traits, such as gene expression profiles \cite{brawand2011evolution},  and have even found use outside biology, such as in internet network tomography \cite{eriksson2010toward,tsang2004network}. Recent work in 
algebraic statistics 
\cite{boege2021reciprocal,sturmfels2020estimating,SUZ20, truell2021maximum} uses algebraic geometry to study parameter inference problems for BMT models. 

The BMT model on a tree has a simple description in terms of the covariance matrices of the densities in the model, stated in  \Cref{def:BMTM} and illustrated in \Cref{fig:brownian_tree}. 
\begin{figure}[ht]
    \centering 
    \begin{subfigure}[c]{.4\textwidth} 
        \centering
        {
\begin{tikzpicture}[scale=1.0,auto=left, thick,->]
    \coordinate (l0) at (0,0);
    \coordinate (v6) at (0,-1);
    \coordinate (v5) at (-1,-2);
    \coordinate (l1) at (-2,-3);
    \coordinate (l2) at (-1,-3);
    \coordinate (l3) at (0,-3);
    \coordinate (l4) at (2,-3);
\tiny
    \draw[green!40!black] (l0)--(v6) node[midway, right, black] {$\theta_0$};
    \draw[green!40!black] (v6)--(v5) node[midway, left, black] {$\theta_5$};

    \foreach \i in {1} {
        \draw[green!40!black] (v5)--(l\i) node[pos=.7, left, black] {$\theta_\i$};
    }
    \foreach \i in {2,3} {
        \draw[green!40!black] (v5)--(l\i) node[pos=.7, right, black] {$\theta_\i$};
    }
    \draw[green!40!black] (v6)--(l4) node[midway, right, black] {$\theta_4$};

    \node[circle, fill=gray!20] at (v5){$5$};
    \node[circle, fill=gray!20] at (v6){$6$};
    \node[circle, fill=red!60!blue!50!white] at (l0){$0$};
    
    \foreach \i in {1,2,3,4} {
        \node[circle, below, fill=red!60!blue!50!white] at (l\i){$\i$};
    }
\end{tikzpicture}} 
    \end{subfigure}%
    \begin{subfigure}[c]{.6\textwidth}
        $\mathcal{L}_\T(\RR) = \left\{ \Sigma=
        \begin{pmatrix} 
            t_1 & t_5 & t_5 & t_6 \\ t_5 
            & t_2 & t_5 & t_6 \\
            t_5 & t_5 & t_3 & t_6 \\
            t_6 & t_6 & t_6 & t_4 
        \end{pmatrix} : t_1, \ldots, t_6 \in \RR \right\}.$
    \end{subfigure}
    \caption{An evolutionary tree $\T$ on  species $1,2,3$ and $4$. A covariance matrix in the associated BMT model must be in the linear space $\mathcal{L}_\T(\RR)$. 
     }
    \label{fig:brownian_tree}
\end{figure}
For simplicity,  label the root of the phylogenetic tree $\T$ by $0$ and the rest of the leaves by $1,\ldots,n$. 
Denote by $\lca(i,j)$ the \textit{least common ancestor} of non-root leaves $i$ and $j$; that is, $\lca(i,j)$ is the first common node on the paths joining $i$ and $j$ to the root $0$. Denote by $\mathrm{Int}(\T)$ the set of all internal nodes of $\T$. 
Whenever $i \neq j$,  $\lca(i,j)$ belongs to $\mathrm{Int}(\T)$. Denote by  $\mathbb{S}^{n}(\mathbb{R})$ the set of $n\times n$ symmetric matrices with real entries and let  $\PD^n(\mathbb R)$ denote the cone of positive definite matrices within $\mathbb{S}^{n}(\mathbb{R}).$ We can now define the BMT model.
\begin{defn}\label{def:BMTM}
Let $\T$ be a phylogenetic tree on $n$ non-root leaves. Consider the linear space of symmetric matrices 
\begin{equation}\label{eq:LSSMdef}
    \mathcal{L}_\T(\mathbb R) := \{\Sigma\in\mathbb{S}^{n}(\mathbb{R})\mid \sigma_{ij}=\sigma_{kl} \text{ if }\lca(i,j)=\lca(k,l)\}. 
\end{equation}
    The \emph{Brownian motion tree model}, $\mathcal{M}_\T$,  specified by $\T$ is the set of all multivariate Gaussian distributions with mean $\mathbf{0}$ and whose covariance matrices lie in the set $\mathcal{L}_\T(\mathbb R)\cap \PD^n(\mathbb R)$. 
\end{defn}
Finding the probability distribution in a fixed model that best fits observed data is a standard problem in statistics and data science. One popular method for inferring such a distribution is maximum likelihood estimation. The \emph{maximum likelihood estimate} (MLE) for given data is the maximizer of the log-likelihood function (see \Cref{sec:likelihood}.1) over the model.
In this paper, we investigate the number of complex critical points of the log-likelihood function over a BMT model. 
This number, known to be invariant under a generic choice of data, is referred to as the maximum likelihood degree (ML-degree) of a model \cite[Chapter~2.1]{drton2008lectures}. Since the MLE, if it exists, is one of these critical points, the ML-degree of a model measures the algebraic complexity of this problem. 

Brownian motion tree models are not exponential families {and are instead examples of \emph{linear Gaussian covariance models}.} As such, their likelihood functions are not typically convex and finding general formulae for their ML-degree is challenging. Sturmfels, Timme, and Zwiernik in  \cite{sturmfels2020estimating}  use numerical algebraic geometry to compute ML-degrees of BMT models for phylogenetic trees with up to $6$ leaves.   

In \cite[Section 4]{SUZ20}, Sturmfels, Uhler, and Zwiernik conjecture a formula for the ML-degree of a star tree model; that is, of a BMT model whose tree has exactly one internal node.
The main result of this paper proves their conjecture. {We phrase this result in terms of a star tree with $n+1$ total leaves; in this setting, one leaf labeled $0$ is the root and there are $n$ leaves representing extant species.}

\begin{theorem}
    \label{thm:main-B}
    The ML-degree of the Brownian motion tree model on a star tree on $n + 1$ leaves is $2^{n+1} - 2(n+1) -1$.
\end{theorem}
    
While computing the ML-degree for trees with multiple internal nodes remains challenging, the next result consolidates the problem to classes of trees with the same unlabeled, unrooted tree~topology. 
    
\begin{theorem}
    \label{thm:main-A}
    Brownian motion tree models on phylogenetic trees with the same unlabeled, unrooted tree topology have the same ML-degree.
\end{theorem}

The proofs of these theorems rely on the toric geometry of  the inverse linear space of $\mathcal{L}_\T$  under the change of variables given by Sturmfels, Uhler, and Zwiernik \cite{SUZ20} (see \Cref{thm:SUZ}), and monomial  parametrization 
provided previously by Boege et al. \cite{boege2021reciprocal} (see \Cref{thm:brownian_toric}). 
This parametrization, called the \emph{path parametrization}, assigns to each edge $e$ in the tree a parameter $\theta_e$ and
allows us to write each concentration matrix $K_\T$ in $\mathcal{L}_\T^{-1}$ as a function of the parameters~$\theta_e$, which we call $K_\T(\theta)$.
Towards this end, we compute the degree of this parametrization in \Cref{thm:degree}, which allows us compute the ML-degree by counting the solutions in the new parameters. 
Though we focus on ML-degrees for star trees in the present work, we envision that \Cref{thm:degree} will be useful for future work towards computing ML-degrees for arbitrary trees.

In order to write the log-likelihood function in terms of the parameters $\theta_e$, we require an expression for the determinant of $K_\T(\theta)$. 
Our result is a weighted analog of the Cayley-Pr\"ufer formula, which is a factorization of the sum of the products of edge variables over spanning trees of an unweighted complete graph.
This yields an explicit formula for $\det(K_\T(\theta))$ for any tree $\T$  (see Theorem \ref{thm:detK-formula}), which is applicable to future work on the likelihood geometry of BMT models.

\subsection{Structure of the paper.} In \Cref{sec:preliminaries}, we reframe Brownian motion tree models in terms of their concentration matrices. We recall the toric representation of this space of concentration matrices by the monomial path map and compute the degree of this map in Theorem \ref{thm:degree}. In \Cref{sec:likelihood}, we introduce the maximum likelihood estimation problem for BMT models and define their ML-degree via the score equations. In order to better understand these score equations, we prove Theorem \ref{thm:detK-formula}, which writes the determinant of an arbitrary concentration matrix in the BMT model in terms of the parameters of the path map.
Section \ref{sec:ResultsParametrization} is devoted to the proof of Theorem \ref{thm:main-A}, which states that the ML-degree is invariant under rerooting.
In \Cref{sec:StarTree}, we prove the ML-degree formula for the BMT model on a star tree.  We end with a discussion of our results and directions for future work.

\section{Toric Geometry of Brownian Motion Tree Models} \label{sec:preliminaries}

\subsection{Monomial parametrizations of concentration matrices}

Recall that the covariance matrices in the Brownian motion tree model specified by a tree $\T$ are exactly the positive definite matrices in $\mathcal{L}_\T$, as described in Definition \ref{def:BMTM}. The \emph{concentration} (or precision) matrices for the BMT model are therefore $\mathcal{L}^{-1}_\T \cap \PD^n(\RR)$, where $\mathcal{L}_\T^{-1}$ is the Zariski closure of all matrices $K = (k_{ij})_{1 \leq i \leq j \leq n}\in \mathbb S^{n}$ such that $K^{-1}\in \mathcal{L}_\T$. {Let $\mathbb{R}[K] = \mathbb{R}[k_{ij} \mid 1 \leq i \leq j \leq n]$ denote the polynomial ring whose variables correspond to entries of $K$. Similarly let $\mathbb{R}(\Sigma)$ denote the fraction field of the polynomial ring with variables $\sigma_{ij}$ for $1 \leq i \leq j \leq n$.} The algebraic variety  $\cL^{-1}_\T$  is the vanishing locus of the kernel of the rational map  
\begin{align} 
    \label{eq:minors}
    \rho_\T: \mathbb{R}[K] \to \mathbb{R}(\Sigma), \ 
    k_{ij} \mapsto \frac{(-1)^{i+j} \det(\Sigma_{ij})}{\det (\Sigma)}, 
\end{align}
where~$\Sigma_{ij}$ is the submatrix of the symmetric matrix~$\Sigma$ with its $i$th row and $j$th column deleted. The ideal $\ker(\rho_\T)$ is referred to as the \emph{vanishing ideal of the Brownian motion tree model}~$\Mcal_\T$. {We note that this is the vanishing ideal of the \emph{inverse} linear space $\mathcal{L}_\T^{-1}$, rather than of the linear space itself.}

\begin{example}\label{ex:KVanishingIdeal}
    Consider the phylogenetic tree $\T$ in Figure \ref{fig:brownian_tree}. It is pictured along with a generic element of its associated linear space $\cL_\T$.
    The inverse linear space $\cL_\T^{-1}$ has vanishing~ideal
    {\footnotesize
     \begin{align*}&\ker(\rho_\T)= \\ &\langle 
        k_{12}k_{14} + k_{14}k_{22} + k_{14}k_{23} - k_{11}k_{24} - k_{12}k_{24} - k_{13}k_{24}, \, k_{13}k_{14} + k_{14}k_{23} + k_{14}k_{33} - k_{11}k_{34} -k_{12}k_{34} - k_{13}k_{34}, \\
        & k_{13}k_{24} + k_{23}k_{24} + k_{24}k_{33} - k_{12}k_{34} - k_{22}k_{34} - k_{23}k_{34}, \, k_{12}k_{23} + k_{12}k_{33} + k_{12}k_{34} - k_{13}k_{22} - k_{13}k_{23} - k_{13}k_{24}, \\
        &  k_{12}k_{13} + k_{12}k_{33} + k_{12}k_{34} - k_{11}k_{23} - k_{13}k_{23} - k_{14}k_{23}, \, k_{12}k_{34} - k_{13}k_{24}, \, k_{12}k_{34} - k_{14}k_{23} \rangle. 
    \end{align*}}
\end{example}

The variety $\cL^{-1}_\T$ is described in detail by Sturmfels, Uhler, and Zwiernik in \cite{SUZ20} as a toric variety.  To exhibit the toric structure of the space of concentration matrices, they show that $\ker(\rho_\T)$ is generated by binomials after a change of coordinates.
 \begin{theorem}[\cite{SUZ20}, 
 Theorem~1.2]\label{thm:SUZ}
     The vanishing ideal for the Brownian motion tree model of the phylogenetic tree $\T$ is toric in variables $p_{ij}$ with $i,j$ distinct elements of $\{0,\dots,n\}$, where 
\begin{equation}
 \label{eq:pcoord} \begin{matrix} 
 p_{ij} &=& - k_{ij} \quad &\quad {\rm for} \,\,i,j > 0, \text{ and } \\
 p_{0i} & = & \sum_{j=1}^n k_{ij} &  {\rm for} \,\,1 \leq i \leq n.
 \end{matrix}
 \end{equation} It is generated by the forms
    $p_{ik} p_{j\ell} - p_{i\ell} p_{jk},$ where $\{i,j\}$ and $\{k,\ell\}$ are cherries in the induced $4$-leaf subtree on any quadruple $i,j,k,\ell\in \{0,\ldots,n\}$.
 \end{theorem}
We note that in Theorem \ref{thm:SUZ}, the subscripts on the variables $p_{ij}$ are unordered so that~$p_{ij} = p_{ji}$. We therefore think of the concentration matrices in the BMT model in the coordinates~$p_{ij}$ for $0 \leq i < j \leq n$ for the rest of the paper. 

Denote by $I_\T$ the ideal $\ker(\rho_\T)$ in the coordinates $p_{ij}$. Since the ideal $I_\T$ is toric, it is the kernel of a monomial map.  Boege et al. \cite{boege2021reciprocal} connect this ideal to the paths in the tree and use these paths to give the monomial parametrization $\varphi_\T$ of this toric variety, called the \emph{path parametrization}: 
 \begin{align} \label{eq:shortest path map}
\varphi_\T: \mathbb{R}[p_{ij} \mid 0\leq i<j\leq n] \to \mathbb{R}[\theta_e \mid e\in E(\T)], \ 
p_{ij} \mapsto \prod_{e \in i \leftrightsquigarrow j} \theta_e, 
\end{align}
where $E(\T)$ is the set of edges of $\T$ and $i \leftrightsquigarrow j$ is the set of edges in the {shortest} path between leaf~$i$ and leaf $j$.

\begin{theorem}[\cite{boege2021reciprocal}, Proposition 3.1]\label{thm:brownian_toric}  The toric vanishing ideal of the Brownian motion tree model on the phylogenetic tree $\T$  in the $p_{ij}$ coordinates defined in Theorem \ref{thm:SUZ} is the kernel of the path map~$\varphi_\T$. \end{theorem}

\begin{example} The Brownian motion tree model for $\T$ in \Cref{fig:brownian_tree} has   toric vanishing ideal 
\begin{align*}
    I_\T = \langle \, &
    p_{01} p_{24} - p_{02} p_{14} ,\,
    p_{01} p_{34} - p_{03} p_{14}  ,\,
    p_{02} p_{34} - p_{03} p_{24} ,\, p_{02} p_{13} - p_{03} p_{12} ,\, \\
    & p_{01} p_{23} - p_{03} p_{12} ,\,
    %
    p_{12} p_{34} - p_{13} p_{24} ,\,
    p_{12} p_{34} - p_{23} p_{14} \,
\rangle.
\end{align*}
Note that the seven generators in the $k_{ij}$ coordinates given in Example \ref{ex:KVanishingIdeal} can be obtained by applying the change of coordinates in Theorem \ref{thm:SUZ} to each of these binomials. This ideal is the kernel of the path map $\varphi_\T$ whose exponent matrix is
\[
A_\T = 
\begin{pNiceMatrix}[first-row,first-col]
    & p_{01} & p_{02} & p_{03} & p_{04} & p_{12} & p_{13} & p_{14} & p_{23} & p_{24} & p_{34} \\ 
        \theta_0 & 1 & 1 & 1 & 1 & 0 & 0 & 0 & 0 & 0 & 0 \\
        \theta_1 & 1 & 0 & 0 & 0 & 1 & 1 & 1 & 0 & 0 & 0\\
        \theta_2 & 0 & 1 & 0 & 0 & 1 & 0 & 0 & 1 & 1 & 0 \\
        \theta_3 & 0 & 0 & 1 & 0 & 0 & 1 & 0 & 1 & 0 & 1 \\
        \theta_4 & 0 & 0 & 0 & 1 & 0 & 0 & 1 & 0 & 1 & 1\\
        \theta_5 & 1 & 1 & 1 & 0 & 0 & 0 & 1 & 0 & 1 & 1
\end{pNiceMatrix}.
\]
\end{example}

The path parametrization proved to be essential in the computation of reciprocal (dual) maximum likelihood degree of BMT models \cite{boege2021reciprocal}. 
It will continue to be instrumental in all of the proofs of the present work.    

\subsection{Degree of the path parametrization} 
Let $\T$ be a tree on $n + 1$ leaves with edge set~$E(\T)$. Let $\mathrm{Lv}(\T) := \{0, \dots, n \}$ denote the leaf set of $\T$. Consider the pullback of the path parametrization:
\begin{align*}
   \varphi_\T^\ast: \mathbb C^{\#E(\T)} \rightarrow   \mathbb C^{\binom{n+1}{2}}, \  (\theta_e)_{e\in E(\T)} \mapsto \left(\prod_{e \in i \leftrightsquigarrow j} \theta_e\right)_{i\neq j\in \mathrm{Lv}(\T)}.
\end{align*}

In the following proposition, we compute the degree of $\varphi_\T^\ast$ by explicitly describing elements in its fibers. 
For each $S\subseteq \mathrm{Int}(\T)$, we define $\epsilon^S$ to be the element of $\{-1,1\}^{\#E(\T)}$ with
\begin{align*}
 \epsilon^S_e=(-1)^{\#(S\cap e)} \ \text{ for each edge } e\in E(\T).
\end{align*}
For any two vectors $u,v$ of the same length $k$, denote by $u\ast v$ their componentwise product; that is, $u \ast v = (u_i v_i)_{i=1}^k$.

  {
\begin{example}
  Let $\T$ be the tree pictured in \Cref{fig:brownian_tree}. In this case, $\mathrm{Int}(\T) = \{5,6\}$, so there are four vectors of the form $\epsilon^S$. These are:
    \begin{align*}
        \epsilon^\emptyset & = (1,1,1,1,1,1), \\
        \epsilon^{\{5\}} & = (1,-1,-1-1,1,-1), \\
        \epsilon^{\{6\}} &= (-1,1,1,1,-1,-1), \text{ and} \\
        \epsilon^{\{5,6\}} &= (-1,-1,-1,-1,-1,1).
    \end{align*}
    We note that for any tree topology, we have that $\epsilon^\emptyset$ is the vector of all ones.
\end{example}}

\begin{theorem}
    \label{thm:degree}
   {{Let $\T$ be a rooted tree.}}  Let $p \in \mathrm{im}(\varphi_\T^\ast)$, with all coordinates non-zero. Let $\hat{\theta}\in \mathbb C^{\#E(\T)}$ such that $\varphi_\T^\ast(\hat{\theta})=p$. {{If  $\T$ has  no degree two nodes,}} then 
    \[(\varphi_\T^\ast)^{-1}(p)=\{\epsilon^S\ast \hat{\theta} \mid S\subseteq \mathrm{Int}(\T) \}.\]
In particular, the degree of $\varphi_\T^\ast$ is $2^{\#\textrm{Int}(\T)}$.
\end{theorem}

{ Before we prove this result, we remark that when $\T$ is a star tree, Theorem \ref{thm:degree} implies that $(\varphi_\T^\ast)^{-1}(p) = \{\hat{\theta}, -\hat{\theta}\}$. 
Indeed, in the star tree case, there is only one internal node, say $u$. When $S = \{u\}$, all edges of $\T$ intersect $S$ in exactly one node, so $\epsilon^{\{u\}}$ is the vector whose entries are all $-1$.}

\iftoggle{EA}{\begin{proof}[Proof outline]
    First we show that for any $\theta, \hat{\theta} \in \left( \varphi_\T^\ast \right)^{-1} (p)$, $\theta_e = \pm \widehat{\theta}_e$ for all $e \in E(\T)$ using rational combinations of the $p_{ij}$'s. Thus, $\theta = \epsilon \ast \widehat{\theta}$ for some $\epsilon \in \{ \pm 1 \}^{\#E(\T)}$. Then, we use induction on the number of internal vertices to show that there exists $S \subseteq \mbox{Int}(\T)$ so that $\epsilon = \epsilon^S$.
\end{proof}}%

{\begin{proof}[Proof of Theorem \ref{thm:degree}]
{{Denote by $F^{p,\hat{\theta}}_\T =(\varphi_{\T}^*)^{-1}(p) \cap \{\epsilon\ast \hat{\theta} \mid \epsilon\in \{\pm 1\}^{\#E(\T)}\}$, the part of the fiber whose coordinates differ from those of $\hat{\theta}$ by a sign.
For ease of readability, we split the proof into two parts. In the first part, we show that for a tree $\T$ with internal nodes of degree at least $3$,   a point in the fiber is of the form $\epsilon \ast \hat{\theta}$ for some $\epsilon\in\{\pm1\}^{\#E(\T)}$; that is,  $(\varphi^*_\T)^{-1}(p)= F^{p,\hat{\theta}}_\T$.  In the second part, we show that for any rooted tree $\T$, a point $\epsilon \ast \hat{\theta}$ in  $F^{p,\hat{\theta}}_\T$ must have  $\epsilon =\epsilon^S$ for some~$S \subseteq \mathrm{Int}(\T)$. }}

\textbf{Part 1}. 
Let $\{u,v\}\in E(\T)$, and the associated coordinate~$\theta_{uv}$ of $\theta\in (\varphi^*_\T)^{-1}(p)$. At least one of the nodes, say $v$, is an internal node. Let $i$ be a leaf such that the path $v\leftrightsquigarrow i$ contains the edge $\{u,v\}$. Let $j,k$ be two other leaves with property that $\{u,v\}$ is not in the paths $v \leftrightsquigarrow j$ or $v\leftrightsquigarrow k$. 
The existence of distinct leaves $i, j, k$ with these properties is guaranteed since $v$ has degree at least three. The condition~$\theta\in (\varphi_\T^\ast)^{-1}(p)$ implies that 
\begin{align}
    p_{ij} =& \label{eq:1} \prod_{e\in i\leftrightsquigarrow v} \theta_e\prod_{e\in v\leftrightsquigarrow j}\theta_e = \prod_{e\in i\leftrightsquigarrow v} \hat{\theta}_e\prod_{e\in v\leftrightsquigarrow j} \hat{\theta}_e\\
    p_{ik} =& \label{eq:2} \prod_{e\in i\leftrightsquigarrow v} \theta_e \prod_{e\in v\leftrightsquigarrow k}\theta_e = \prod_{e\in i\leftrightsquigarrow v} \hat{\theta}_e\prod_{e\in v\leftrightsquigarrow k} \hat{\theta}_e\\
    p_{jk} =& \label{eq:3}\prod_{e\in j\leftrightsquigarrow v} \theta_e\prod_{e\in v\leftrightsquigarrow k}\theta_e = \prod_{e\in j\leftrightsquigarrow v} \hat{\theta}_e\prod_{e\in v\leftrightsquigarrow k} \hat{\theta}_e.
\end{align}
Since each $p_{ij}\neq 0$, each of the products in the above equations is nonzero.  Solving for~$\prod\limits_{e\in v\leftrightsquigarrow j}\theta_e $ in~(\ref{eq:1}), for $\prod\limits_{e\in v\leftrightsquigarrow k}\theta_e $ in (\ref{eq:2}), and substituting these expressions in (\ref{eq:3}), gives 
\begin{align}
    & \label{eqn:squares} %
    \left( \prod\limits_{e\in i\leftrightsquigarrow v} \theta_e \right)^2 = \left( \prod\limits_{e \in i \leftrightsquigarrow v} \hat{\theta}_e \right)^2.
\end{align}
Therefore, the product of $\theta_e$'s over the path $v \leftrightsquigarrow i$ from an internal vertex $v$ to any leaf $i$ is equal to the product of $\hat{\theta}_e$'s over $v \leftrightsquigarrow i$ up to a sign. Factoring out the terms corresponding to the edge $\{u,v\}$ gives
\begin{align} \label{eq:4}
    \theta_{uv}\prod_{e\in i\leftrightsquigarrow u}\theta_e = \pm  {\hat{\theta}}_{uv}\prod_{e\in i\leftrightsquigarrow u} \hat{\theta}_e.  
\end{align}
If $u=i$ is a leaf, then we immediately obtain $\theta_{uv} = \epsilon_{uv}{\hat{\theta}}_{uv}$ for some $\epsilon_{uv}\in \{-1,1\}$ and we are done. If $u$ is an internal node, then analogously to \Cref{eqn:squares} we have
\begin{align} \label{eq:5}
    \prod_{e\in i\leftrightsquigarrow u}\theta_e= \pm  \prod_{e\in i\leftrightsquigarrow u} \hat{\theta}_e, 
\end{align}
and \Cref{eq:4} and \Cref{eq:5} imply $\theta_{uv}=\epsilon_{uv}{\hat{\theta}}_{uv}$ for some $\epsilon_{uv}\in \{-1,1\}$, as desired.\\


\textbf{Part 2.} {{We will use induction on the number of internal nodes to show that 
\begin{align}\label{eq:part2}
   F^{p,\hat{\theta}}_\T=\{\epsilon^S \ast \hat{\theta}\mid S\subseteq \mathrm{Int}(\T)\}.
\end{align}}}
The following observation will be useful:
\begin{align}\label{eq:claim}
  \text{``if  $\epsilon\ast \hat{\theta}\in F^{p,\hat{\theta}}_\T$, then $\epsilon_{iu}$ is fixed for all leaves $i$ {{that are adjacent to}} $u$".}  
\end{align}
Indeed, let $i,j$ be leaves that are adjacent to $u$. The path between $i$ and $j$ in $\T$ consists of the edges $\{i,u\}$ and $\{j,u\}$. Since $\epsilon\ast \hat{\theta}$ and $\hat\theta$ belong to the same fiber of $\varphi^\ast_\T$, by definition of the path map we have $(\epsilon_{iu}\hat\theta_{iu})(\epsilon_{ju}\hat\theta_{ju}) = \hat\theta_{iu}\hat\theta_{ju}.$
Hence, $\epsilon_{iu} = \epsilon_{ju}$.

{{Now we prove \Cref{eq:part2} by induction on the number of internal nodes. First, if $\T$ has no internal nodes, then $ F^{p,\hat{\theta}}_\T=\{\epsilon^\emptyset\ast \hat\theta\}$. Let $\T$ have one internal node. Call this node~$u$. Let  $\epsilon\ast\hat\theta\in F^{p,\hat{\theta}}_\T$ and suppose $\epsilon_{iu}=-1$ for some leaf $i$. For any  leaf $j$ in $\T$, by Observation (\ref{eq:claim}), one has $\epsilon_{iu} = \epsilon_{ju}=-1$. 
 So, $\epsilon=\epsilon^{\{u\}}$ and $F^{p,\hat{\theta}}_\T=\{\epsilon^{\emptyset}\ast\hat{\theta}, \epsilon^{\{u\}}\ast \hat{\theta}\}$. }}

Suppose that Equation (\ref{eq:part2}) holds for any tree  with less than $m>1$ internal nodes, and let $\T$ have $m$ internal nodes. Select an internal node $u$ in $\T$ whose children are all leaves, which we label $1,2,\ldots, k$ without loss of generality. Note that since $\T$ is not a star tree, any such $u$ must have a parent $v$ that is an internal node.

Let $\T'$ be the tree obtained by removing from $\T$ the leaves $1,\ldots,k$ along with all edges containing $u$. The tree $\T'$ has $m-1$ internal nodes and $\mathrm{Lv}(\T')\subset \mathrm{Lv}(\T)$. Take $p'$ to be the projection of  $p$ in $\CC^{\binom{n+1-k}{2}}$ {{onto the coordinates indexed by pairs in $\Lv(\T')$, and let $\hat{\theta}'=\hat{\theta}|_{E(\T')}$ be the projection of $\hat{\theta}$ onto the edges of~$\T'$. By the induction hypothesis, we have that 
\begin{align}\label{eq:part2 T'}
   F^{p',\hat{\theta}'}_{\T'}=\{\epsilon^{S'} \ast \hat{\theta}'\mid S'\subseteq \mathrm{Int}(\T)\setminus\{u\}\}.
\end{align}
 For any $\epsilon \ast \hat\theta\in F^{p,\hat{\theta}}_\T$, the projection ~$(\epsilon \ast \hat\theta)|_{E(\T')}$  is in $F^{p',\hat{\theta}'}_{\T'}$.  In other words, every point in  $F^{p,\hat{\theta}}_\T$ can be obtained by extending a point in~$F^{p',\hat{\theta}'}_{\T'}$. In the following paragraph, we show that each point $\epsilon^{S'}\ast \hat{\theta}'\in F^{p',\hat{\theta}'}_{\T'}$ for $S'\subseteq \mathrm{Int}(\T)\setminus \{u\}$ can be extended to precisely the two points $\epsilon^{S'}\ast \hat{\theta}$ and $\epsilon^{S'\cup \{u\}}\ast \hat{\theta}$ in $ F^{p,\hat{\theta}}_{\T}$.  The union of all such points is all $F^{p,\hat{\theta}}_{\T}$. 

Fix some $S'\subseteq \mathrm{Int}(\T)\setminus\{u\} $. 
Then $\epsilon^{S'} \ast \hat{\theta}'\in F^{p',\hat{\theta}'}_{\T'}$. We extend it to a point $\theta$ in  $\mathbb C^{\#E(\T)}$ as follows.
From Observation~\ref{eq:claim}, any point $\epsilon\ast \hat \theta$ in the fiber of $p$ has fixed~$\epsilon_{iu}$ for $i=1,\ldots,k$. 
We call this value $\epsilon_u$ and fix $\epsilon_u \in \{ \pm 1\}$.
The value $\theta_{uv}$ is the only entry of $\theta$ which has not yet been specified.
For $\theta$ to be in $F^{p,\hat{\theta}}_{\T}$, we must have \[\theta_{uv}=\left(\epsilon_u  \prod\limits_{e\in 0\leftrightsquigarrow v}\epsilon^{S'}_e \right)  \hat\theta_{uv}\] since the product of the entries of $\theta$ along the path from $0$ to any leaf $i=1,\dots,k$ must be equal to the same product of the entries of $\hat{\theta}$. 
So, $\theta$ is uniquely determined with:}}
\begin{align}\label{eq:thetas induction}
\theta_e=
\begin{cases}
\epsilon^{S'}_e \cdot \hat{\theta}_e & \text{ for } e\in E(\T'),\\
\epsilon_u \cdot \hat\theta_{e} & \text{ for } e=\{i,u\}, \ i=1,\ldots,k, \text{ and}\\
\left(\epsilon_u  \prod\limits_{e\in 0\leftrightsquigarrow v}\epsilon^{S'}_e \right)  \hat\theta_{e} & \text{ for } e=\{u,v\}.
 \end{cases}
\end{align}
By construction, $\theta=\epsilon^S\ast \hat \theta$, where $S= S'$ when $\epsilon_u=1$ and $S= S'\cup \{u\}$ when $\epsilon_u=-1$. 
{{Parts 1 and 2 together prove that when all internal nodes of $\T$ have degree at least 3, then $(\varphi^{*}_{\T})^{-1}(p)=F^{p,\hat{\theta}}_{\T}=\{\epsilon^S \ast \hat{\theta}\mid S\subseteq \mathrm{Int}(\T)\}$, as desired.}}
\end{proof}

From \Cref{thm:degree} we deduce \Cref{cor:F(k) F(theta)} which is useful for solving polynomial systems arising from Brownian motion tree models such as the likelihood equations and other optimization problems. 
 
\begin{cor}
    \label{cor:F(k) F(theta)}
    Let 
    $F = \{f_1(K),\ldots, f_\ell(K)\}$ be a polynomial system in  the variables $p_{ij}$. 
    Let $F\circ \varphi_\T=\{f_1(K(\theta)),\ldots, f_\ell(K(\theta))\}$ be these polynomials written in variables $\theta_e$. Then, 
    \begin{enumerate}
        \item [a. ] the solutions $K\in \mathcal{L}_\T^{-1}$ to system $F$  are precisely   $K(\theta)$ for $\theta$ a solution of $F\circ \varphi_\T$,
        \item [b. ] \( \deg \big( \langle F \rangle + I_\T \big) = \dfrac{\deg \big(\langle F \circ \varphi_\T\rangle \big)}{2^{\# \mathrm{Int}(\T)}} \).
    \end{enumerate}
\end{cor}

\section{Score Equations}
\label{sec:likelihood}

\subsection{Maximum likelihood estimation in Brownian motion tree models}

Maximum likelihood estimation is a method for inferring the distribution in a statistical model that best explains a data set. 
Let $\mathbf{u}_1,\ldots,\mathbf{u}_m\in \mathbb R^n$ be independent, identically distributed data which we assume are sampled from a distribution in the BMT model on a tree $\T$. The observed data has sample covariance matrix
\begin{equation}
S := \frac{1}{m}\sum^m_{j=1}\mathbf{u}_i\mathbf{u}_i^T\in \PSD^n.
\label{eqn:SCov}
\end{equation}

A \emph{maximum likelihood estimate} (MLE) for this data in the BMT model $\mathcal M_\T$ is a 
concentration matrix $\hat{K} \in \mathcal{L}^{-1}_\T \cap \, \PD^n(\RR)$ that maximizes the value of the density function for the normal distribution $\mathcal{N}(\mathbf{0}, K^{-1})$ on this data, if such a maximizer exists.
Equivalently, $\hat{K}$ is a global maximizer of 
\begin{align}
        \label{eqn: log-likelihood K}
    \ell(K|S)  := \log\det(K)-\trace(SK).
\end{align}
We note that, as written, the expression $\ell(K|S)$ is not exactly the logarithm of the likelihood function. However, they only differ by constant addition and multiplication and hence have the same critical points.  So we slightly abuse terminology and refer to $\ell(K|S)$ as the \emph{log-likelihood function}. While the function $\ell(K|S)$ is convex over the positive definite cone, it is not convex when $K$ is restricted to the BMT model \cite{sturmfels2020estimating}. We refer the reader to \cite[Chapter~7]{sullivant2018algebraic} for background on algebraic geometry and maximum likelihood estimation. Section 2 of \cite{coons2020maximum} also thoroughly introduces maximum likelihood estimation, specifically in linear covariance models.

The MLE is a critical point of the log-likelihood function. Hence, we may compute it by finding the common zeros of the partial derivatives of $\ell(K|S)$ and computing the likelihood at each critical point. Thus, the number of critical points of $\ell(K|S)$, called the \emph{maximum likelihood degree} (ML-degree), measures the algebraic complexity of computing the MLE. We now define the ML-degree more precisely.

\begin{defn}\label{def:MLD}
    The \emph{maximum likelihood degree} of the BMT model, denoted~$\mld(\mathcal{M}_\T)$, is the number of complex critical points $\ell(K|S)$ over $\mathcal{M}_\T$, counted with multiplicity, for a generic sample covariance matrix $S$.
\end{defn}

In order to compute the ML-degree of a Gaussian model, we begin by writing the log-likelihood $\ell(K|S)$ in terms of the parameters $\theta_e$ of the path map. Since we are interested in the critical points of the log-likelihood, we take the partial derivatives of $\ell(K|S)$ with respect to each $\theta_e$ and set these equal to zero. These partial derivatives are called the \emph{score equations}.
In the case of a linear Gaussian covariance model, they are rational functions. In fact, they are of the form
\[
\frac{\partial \ell}{\partial \theta} = \frac{1}{\det(K)} \frac{\partial}{\partial \theta}(\det(K)) - \frac{\partial}{\partial \theta}(\trace(SK)).
\]
We can compute the vanishing locus of the rational score equations by finding the variety of their numerators and removing the variety of the product of their denominators.
In the Gaussian case, note that since $\det(K)$ and $\trace(SK)$ are both polynomials in the $\theta_e$ parameters, the only denominator that appears in any score equation is $\det(K)$. Hence, removing the vanishing locus of the denominators simply corresponds to removing any solutions for which the resulting concentration matrix would be singular. Next, we count the critical points in this variety in the $\theta_e$ coordinates, including their multiplicities.  Finally, Theorem \ref{thm:degree} allows us to divide the number of solutions in the $\theta_e$ parameters by $2^{\#\textrm{Int}(T)}$ to obtain the ML-degree of the BMT model.

\begin{example}
Consider the tree $\T$ in Figure \ref{fig:brownian_tree}. Via the change of coordinates in Theorem \ref{thm:SUZ} and the path map, the concentration matrices in the BMT model on $\T$ are of the following form in the $\theta_e$ parameters:
    \[ \tiny \begin{pmatrix}
        \theta_1 \left( \theta_5 \left( \theta_0 + \theta_4 \right) + \theta_2 + \theta_3 \right) & -\theta_1 \theta_2 & -\theta_1 \theta_3 & -\theta_5 \theta_1 \theta_4 \\
        -\theta_1 \theta_2 & \theta_2 \left( \theta_5 \left( \theta_0 + \theta_4 \right) + \theta_1 + \theta_3 \right) & -\theta_2 \theta_3 & -\theta_5 \theta_2 \theta_4 \\
        -\theta_1 \theta_3 & -\theta_2 \theta_3 & \theta_3 \left( \theta_5 \left( \theta_0 + \theta_4 \right) + \theta_1 + \theta_2 \right) & -\theta_5 \theta_3 \theta_4 \\
        -\theta_5 \theta_1 \theta_4 & -\theta_5 \theta_2 \theta_4 & -\theta_5 \theta_3 \theta_4 & \theta_4 \left( \theta_0 + \theta_5 \left( \theta_1 + \theta_2 + \theta_3 \right) \right)
    \end{pmatrix}. \]
Consider a generic sample covariance matrix $S = (s_{ij})_{1 \leq i ,j \leq 4}$.
The log-likelihood function for~$S$ in this BMT model is
\begin{align*}
    \ell(K|S) &= \log \det(K) - \trace(SK) \\
    &= \log(\theta_0\theta_1\theta_2\theta_3\theta_4\theta_5(\theta_1\theta_5 + \theta_2\theta_5 + \theta_0 + \theta_4)(\theta_0 \theta_5 + \theta_4 \theta_5 + \theta_1 + \theta_2 + \theta_3)^2) \\ & \quad
     - s_{11} \theta_1 \left( \theta_5 \left( \theta_0 + \theta_4 \right) + \theta_2 + \theta_3 \right)  - s_{22}\theta_2 \left( \theta_5 \left( \theta_0 + \theta_4 \right) + \theta_1 + \theta_3 \right) \\ & \quad - s_{33}\theta_3\left( \theta_5 \left( \theta_0 + \theta_4 \right) + \theta_1 + \theta_2 \right) - s_{44}\theta_4 \left( \theta_0 + \theta_5 \left( \theta_1 + \theta_2 + \theta_3 \right) \right) \\ 
     & \quad + 2s_{12} \theta_1\theta_2 + 2 s_{13} \theta_1 \theta_3 + 2s_{14} \theta_1 \theta_4 \theta_5 + 2s_{23} \theta_2 \theta_3 + 2 s_{24} \theta_2 \theta_4 \theta_5 + 2s_{34} \theta_3 \theta_4 \theta_5.
\end{align*}
{
The score equations are the partial derivatives of the function above. For example, the first score equation is:
{\small $$\frac{\partial \ell}{\partial \theta_0} = \frac{1}{\theta_0} + \frac{1}{\theta_1 \theta_5 + \theta_2 \theta_5 + \theta_0 + \theta_4} + \frac{2 \theta_5}{\theta_0 \theta_5 + \theta_4 \theta_5 + \theta_1 + \theta_2 + \theta_3} - s_{11} \theta_1 \theta_5 - s_{22} \theta_2 \theta_5 - s_{33} \theta_3 \theta_5 - s_{44} \theta_4.$$}}
To find the ML-degree, we need to set the system of score equations equal to zero and solve. We did this using the \texttt{Julia} software package \texttt{LinearCovarianceModels.jl}, which makes use of homotopy continuation and has functionality specifically designed to compute ML-degrees of Brownian motion tree models \cite{HomotopyContinuation,sturmfels2020estimating}. We refer the reader to the documentation of \texttt{LinearCovarianceModels.jl} for precise instructions on how to use this package to construct BMT models and calculate their ML-degrees. 
From these computations, we see that the system has $44$ solutions in $(\CC^\ast)^6$. By Theorem \ref{thm:degree}, the degree of the path map in this case is $4$,  since $\T$ has two interior nodes. Hence, dividing by $4$ gives that~$\mld(\mathcal{M}_\T)=11.$
\end{example}

\subsection{A generalization of the Cayley-Pr\"ufer Theorem}
\label{subsec:cayley-prufer}

The classical Cayley-Pr\"ufer Theorem provides  an enumeration of the spanning trees of a complete graph in factored~form
\begin{theorem}[Cayley-Pr\"ufer Theorem, \cite{cayley1878theorem}]
    \label{thm: classical cayley prufer}
    Let $\Kcal_n$ be the complete graph on $n$ vertices. Then,
    \begin{equation}
        \label{eqn:cp}
        \sum_{\substack{\Gamma \subseteq \Kcal_n \\ \textrm{spanning}\\\textrm{tree}}} \prod_{v \in V(\Kcal_n)} x_v^{\deg_\Gamma(v)} = x_1 \cdots x_n \left( x_1 + \cdots + x_n \right)^{n - 2},
    \end{equation}
    where $V(\Kcal_n) = [n]$ is the vertex set of $\Kcal_n$ and $\deg_\Gamma(v)$ is the number of edges adjacent to $v$ in the tree~$\Gamma$. 
\end{theorem}

The goal of this section is to prove \Cref{thm:detK-formula}, which factorizes $\det (K_\T(\theta))$ and specializes to the classical Cayley-Pr\"ufer Theorem when $\T$ is a star tree.
We begin by recalling Kirchoff's Matrix-Tree Theorem.
Let $\Gcal$ be a weighted graph with vertex set $[n]:=\{1,\ldots,n\}$. 
Let $w_{ij}$ be the weight for the edge $\{ i, j \} \in E(\Gcal)$. 
One naturally extends $w$ to all pairs of vertices in $\Gcal$ by setting $w_{ij} = 0$ when $\{ i, j \}$ is not an edge of $\Gcal$.
The \emph{weighted Laplacian of $\Gcal$}, denoted $L_\Gcal$, is an $n \times n$ matrix which encodes the weights of $\Gcal$ as follows:
\[ (L_\Gcal)_{ij} = 
    \begin{cases}
        \displaystyle\sum_{\substack{k = 1 \\ k \neq i}}^n w_{ik} & \text{ if } i = j \text{ and} \\
        - w_{ij} & \text{ if } i \neq j.
    \end{cases}
\]

Kirchoff's Matrix-Tree Theorem \cite{kirchoff}, applied to a complete graph {with the edge weights $w_{ik} = x_i x_k$,} states that the left hand-side sum in \Cref{eqn:cp} is the determinant of any principal submatrix of the Laplacian of $\Kcal_n$.
Our concentration matrix, $K_\T(\theta)$, is a principal submatrix of the Laplacian of a  weighted complete graph, $\Kcal^\T$ with weights determined by the paths in tree $\T$.  

\begin{defn}
    \label{defn:Kcal-T}
    Let $\T$ be a phylogenetic tree on $n + 1$ leaves. Define $\Kcal^\T$ to be the weighted complete graph on $n + 1$ vertices, where the weight of an edge $\{ i, j \}$ is  $\varphi_\T(p_{ij}) = 
    \prod\limits_{e \in i \leftrightsquigarrow j} \theta_e$.
\end{defn}

\begin{figure}[h]
        \centering
        \begin{tikzpicture}[scale=2.5]
  \tikzstyle{vertex}=[circle, draw, fill=black!25, minimum size=20pt, inner sep=0pt]
  
    \foreach \i in {0,...,4} {
        \coordinate (v\i) at ({90 + 360/5 * (\i)}:1);
    }
  
    \foreach \i in {0,...,4} {
        \foreach \j in {0,...,4} {
            \draw (v\i) -- (v\j);
        }
    }

    \draw (v0)--(v4) node[midway, above right] {$\varphi_\T(p_{04}) = \theta_0 \theta_4$};
    \draw (v0)--(v1) node[midway, above left] {$\varphi_\T(p_{01}) = \theta_0 \theta_1 \theta_5$};
    \draw (v1)--(v2) node[midway, below left] {$\varphi_\T(p_{12}) = \theta_1 \theta_2$};

    \foreach \i in {0,...,4} {
        \node[vertex] at (v\i){\i};
    }
\end{tikzpicture}
        \caption{The weighted complete graph, $\Kcal_{5}^\T$, for the tree from \Cref{fig:brownian_tree}.}
        \label{fig:Kcal}
    \end{figure}

The Laplacian of $\Kcal^\T$ is the $n \times n$ matrix with entries
\[
    (L_{\Kcal^\T})_{0 \leq i, j \leq n} = 
    \begin{cases}
        \displaystyle \sum_{\substack{k = 0 \\ k \neq i}}^n \prod_{e \in i \leftrightsquigarrow k} \theta_e & \text{if } i = j \text{ and} \\
        \displaystyle - \prod_{e \in i \leftrightsquigarrow j} \theta_e & \text{if } i \neq j.
    \end{cases}
\]
Note that $K_\T(\theta)$ is the submatrix of $L_{\Kcal^\T}$ with row and column corresponding to vertex 0 removed.
By the Matrix-Tree Theorem, $\det (K_\T(\theta))$ enumerates the weighted spanning trees of $\Kcal^\T$. The following theorem provides a factorization for this enumeration.

\begin{theorem} 
    \label{thm:detK-formula} Let $\T$ be any phylogenetic tree. The determinant of the concentration matrix $K_\T(\theta)$ in the BMT model on $\T$ is 
    \begin{equation}
        \label{eqn:detK-formula}
        \det (K_\T(\theta)) = \left( \prod_{e \in E(\T)} \theta_e \right) \prod_{v \in \mathrm{Int}(\T)} \left( \sum_{\ell \in \mathrm{Lv}(\T)} \prod_{e \in v \leftrightsquigarrow \ell} \theta_e \right)^{\text{deg}_{\T}(v) - 2}.
    \end{equation}
\end{theorem}

\begin{proof}
    Note that since $K_\T$ is a principal submatrix of the Laplacian of a complete graph, its determinant is described by the matrix-tree theorem:
    \begin{equation*}
        \label{defn:weighted spanning tree enum}
        \det (K_\T(\theta)) = \sum_{\substack{\Gamma \subset \Kcal^\T \\ \textrm{spanning} \\ \textrm{tree}}} \prod_{\{i,j\} \in E(\Gamma)} \varphi_\T(p_{ij}).
    \end{equation*}

    Note that $\det (K_\T(\theta))$ is homogeneous of degree $2n$ in the variables $\theta_e$ where $e$ contains a leaf of $\T$.
    We call these variables \textit{leaf variables} and introduce the following notation.
    For each leaf $i$, let $v_i$ be the vertex of $\T$ that shares an edge with $i$. 
    We use $\theta_i$ to denote $\theta_{\{ v_i, i\}}$.
    
    Now observe that $\theta_e$ divides $\det (K_\T (\theta))$ for all edges~$e \in E(\T)$. 
    {Indeed, for every edge $e$ in the tree $\T$, any spanning tree of $\Kcal^\T$ must contain at least one edge $\{i, j\}$, where the leaves $i$ and $j$ are in different connected components of $\T \setminus \{ e \}$ since otherwise $\Kcal^\T$ would be disconnected}.
    Thus, $\theta_e$ divides each term of the sum in~$\det (K_\T (\theta))$. 
    Define
  \begin{align*}
    D(\T) &:= \prod_{v \in \mathrm{Int}(\T)} \left( \sum_{\ell \in \mathrm{Lv}(\T)} \prod_{e \in v \leftrightsquigarrow \ell} \theta_e \right)^{\text{deg}_\T(v) - 2} \text{and} \\
    \Delta(\T) &:= \frac{\det (K_\T (\theta))}{\prod\limits_{e \in E(\T)} \theta_e} - D(\T).
  \end{align*}

    To prove the theorem, we will show that $\Delta(\T) = 0$ by induction on the number of leaves.
    If $\T$ has two leaves, then it is easy to check that $\Delta(T) = 0$.
    {For $\T$ a tree on $n+1$ leaves, $\Delta(\T)$ is homogeneous of degree $n - 1$ in the leaf variables.
    To see this, note that the summands in the definition of $\Delta(\T)$ are homogeneous of degree $n-1$ in the leaf variables.
    Indeed, the degree of the first summand in the leaf variables is $2n - (n+1) = n - 1$.
    For the second summand, note that $\sum_{v \in \In(\T)} \deg_\T(v) = \# E(\T) - \#\Lv(T)$, since the sum counts internal edges twice and leaf-adjacent edges once.
    Furthermore, since $\T$ is a tree we have $\# \In(\T) = \# E(\T) - \# \Lv(\T) + 1$.
    It then follows that
    \[ \sum_{v \in \In(\T)} (\deg_\T(v) - 2) = \# E(\T) - \#\Lv(T) - 2 (\# E(\T) - \# \Lv(\T) + 1) = n - 1, \]
    and so the second summand also has degree $n - 1$ in the leaf variables.}
     
{    On the other hand, we show below that $\Delta(\T)$ is divisible by each of the $n + 1$ leaf variables.
    Assuming this, and if $\Delta(\T) \neq 0$, then $\Delta(\T)$ contains at least one term $m \cdot \theta_0^{a_0} \theta_1^{a_1} \cdots \theta_n^{a_n}$, where  $m$ is a monomial in the non-leaf variables and $a_i > 0$ for $i = 0, 1, \ldots, n$.
    Hence, the degree of the monomial in the leaf variables is at least $n + 1$.
    This contradicts the fact that~$\Delta(\T)$ is homogeneous of degree $n - 1$ in the leaf variables, and so $\Delta(\T)$ must be identically zero.
    Thus, it suffices to prove the claim that $\Delta(\T)$ is divisible by each of the $n + 1$ leaf variables.}
    
    Fix a leaf $i \in \mathrm{Lv}(\T)$. 
    We show that $\theta_i$ divides~$\Delta(\T)$, or equivalently, $\Delta(\T)\vert_{\theta_i = 0} = 0$. We claim that
    \begin{equation}
        \label{eqn: delta_T}
        \Delta(\T) \bigg\rvert_{\theta_i = 0} = \Delta(\T \setminus \{i\}) \left( \sum_{j \in \mathrm{Lv}(\T) \setminus \{i\}} \frac{\varphi_\T(p_{ij})}{\theta_i} \right).
    \end{equation}
    The tree $\T \setminus \{ i \}$ has fewer leaves than $\T$, so by induction, $\Delta(\T \setminus \{ i \}) = 0$. Thus, if \Cref{eqn: delta_T} holds, $\theta_i$ divides $\Delta(\T)$. To prove \Cref{eqn: delta_T}, we will show the following:
    \begin{align}
        \label{eqn: F induction}
        \frac{\det (K_\T)}{\prod\limits_{e \in E(\T)} \theta_e} \Bigg\rvert_{\theta_i = 0} &= \frac{\det (K_{\T \setminus \{i\}})}{\prod\limits_{e \in E(\T) \setminus \{i\}} \theta_e} \cdot \left( \sum_{j \in \mathrm{Lv}(\T) \setminus \{i\}} \frac{\varphi_\T(p_{ij})}{\theta_i} \right) \text{ and} \\
        \label{eqn: D induction}
        D(\T) \bigg\rvert_{\theta_i = 0} &= D(\T \setminus \{i\}) \cdot \left( \sum_{j \in \mathrm{Lv}(\T) \setminus \{i\}} \frac{\varphi_\T(p_{ij})}{\theta_i} \right).
    \end{align}

    On the left-hand side of \Cref{eqn: F induction}, the term corresponding to the spanning tree $\Gamma$ has degree $\deg_\Gamma(i) - 1$ in $\theta_i$.
    Therefore, the terms on the left-hand side that remain after setting $\theta_i = 0$ are those where $i$ is a leaf in $\Gamma$. 
    The factorization on the right says that any such spanning tree is obtained by first finding a spanning tree on $\Kcal^\T \setminus \{i\}$, and then adding $i$ as a leaf. This proves \Cref{eqn: F induction}.

    To prove \Cref{eqn: D induction}, let $v_i$ be the internal vertex of $\T$ adjacent to $i$. Then
    $$D(\T) =  \left( \theta_i + \sum_{\ell \in \mathrm{Lv}(\T) \setminus \{ i \}} \frac{\varphi_\T(p_{i \ell})}{\theta_i} \right)^{\text{deg}(v_i) - 2} \cdot \prod_{v \in \mathrm{Int}(\T) \setminus \{v_i\}} \left( \sum_{\ell \in \mathrm{Lv}(\T)} \prod_{e \in v \leftrightsquigarrow \ell} \theta_e \right)^{\text{deg}(v) - 2}.$$
    When $\theta_i = 0$, the leftmost factor is 
    \[\sum_{j \in \mathrm{Lv}(\T) \setminus \{i\}} \frac{\varphi_\T(p_{ij})}{\theta_i},\] and the product of the rest of the factors is $D({\T \setminus \{i\}})$, which proves \Cref{eqn: D induction} and concludes the proof.
\end{proof}

We note that when $\T$ is a star tree, \Cref{thm:detK-formula} recovers the classical Cayley-Pr\"ufer formula (\Cref{thm: classical cayley prufer}).

\iftoggle{EA}{}{\begin{cor}
    \label{prop:star-detK-formula}
    Let $\T$ be the star tree on $n+1$ leaves, with each edge is weighted by $\theta_i$, where $i$ is the adjacent leaf. Then,
    \begin{equation}
    \label{eqn:star-detK-formula}
        \det (K_\T(\theta)) = \prod_{i=0}^n \theta_i \left( \sum_{i=0}^n \theta_i \right)^{n-1}.
    \end{equation}
\end{cor}}

\section{Equivalence of phylogenetic trees up to re-rooting}
\label{sec:ResultsParametrization}

In this section, we show that the ML-degree of a Brownian motion tree model depends only on the (unlabeled) unrooted tree topology. In particular, we show that the ML-degree does not depend on which of the $n+1$ leaves is chosen to be the root. Moreover, we show that if $\T$ and $\T^\prime$ are two trees with the same unrooted tree topology, the MLE of $\T^\prime$ can be easily obtained from the MLE of $\T$.

Let $\T$ be a phylogenetic tree on leaves $\mathrm{Lv}(\T)=\{0,\ldots,n\}$ with $0$ as its root. Let $r \in [n]$ and let $\T'$ be the rooted phylogenetic tree with the same unrooted topology as $\T$ obtained by rerooting at $r$. The non-root leaves of $\T'$ are then $\mathrm{Lv}(\T)\setminus \{r\}$. We consider them in the order~$1,\dots,r-1,0,r+1,\dots,n$. With this order, an arbitrary element $K'$ of $\cL_{\T'}^{-1}$ has entries 
\begin{align}\label{eqn:K'}
    k'_{ij}=\begin{cases}  \sum\limits_{t = 1}^n p_{0t}  & \text{ for } i= j=r,\\
          -p_{0j} & \text{ for } i=r \text{ and } j\neq r,\\
    \sum\limits_{\substack{t = 0 \\ t \neq i}}^n p_{it}  & \text{ for } i= j\neq r, \text{ and}\\
       -p_{ij} & \text{ for  } i\neq j \text{ and } i,j\neq r,
    \end{cases}
\end{align}
where  $p_{ij}$ are as in Equation (\ref{eq:pcoord}). This gives an invertible linear transformation between the varieties $\cL^{-1}_{\T}$ and~$\cL^{-1}_{\T'}$.

Given a symmetric matrix $S$, construct the symmetric matrix  $S'$  by  applying the following invertible linear transformation to the entries of $S$:
\begin{align}\label{eq:S'}
\begin{split}
  s'_{rr}&=s_{rr},\\
      s'_{rj} &=s_{rr}-s_{rj} \text{ for } j\neq r,\\
    s'_{ii}&=s_{rr}+s_{ii}-2s_{ri} \text{ for } i\neq r,\\
    s'_{ij}&= s_{rr}-s_{ri}-s_{rj}+s_{ij}  \text{ for } i\neq j \text{ and } i,j\neq r.
    \end{split}
\end{align}

This linear transformation is visibly invertible since it writes the vectorization of $S'$ as an upper-triangular matrix times the vectorization of $S$.
Let $\mathrm{mle}(\mathcal{M}_\T,S)$ denote the MLE for a sample covariance matrix $S$ in $\mathcal{M}_\T$. We consider $\mathrm{mle}(\mathcal{M}_\T,S)$ to be written in the coordinates $(p_{ij})_{0 \leq i < j \leq n}$.

\begin{reptheorem}{thm:main-A}
    \label{thm:change_of_root} 
    Let $\T$ and $\T'$ be phylogenetic trees with the same unlabeled, unrooted tree topology as $\T$. Then
    \begin{enumerate}
        \item [(a)] $\mld(\mathcal{M}_{\T})=\mld(\mathcal{M}_{\T'})$ and
        \item [(b)] $\mathrm{mle}(\mathcal{M}_\T,S)=\mathrm{mle}(\mathcal{M}_{\T'},S')$ for $S$ and $S'$ as in \Cref{eq:S'}, if both MLEs exist. 
    \end{enumerate}
\end{reptheorem}

\iftoggle{EA}{\begin{proof}[Proof outline]
    We prove $\ell_\T(K|S) = \ell_{\T'}(K^\prime|S^\prime)$ in two steps. First, \Cref{thm:detK-formula} gives $\det K = \det K^\prime$. Second, we show that $\mbox{tr} (SK) = \mbox{tr}(S^\prime K^\prime)$. We conclude that we can recover the critical points of the log-likelihood function for $\T$ from the log-likelihood function of $\T^\prime$, and vice-versa. Finally, we show that if $K$ is positive definite, then $K'$ is as well.
\end{proof}}%
{\begin{proof}
    Assume without loss of generality that $\T$ has root $0$  and non-root leaf set $[n]$, and $\T'$ has root $r \in \mathrm{Lv}(\T)$ and the same unlabeled, unrooted topology as $\T$. Let $K \in \cL_{\T}^{-1}$ in the $p_{ij}$ coordinates as in Equation (\ref{eq:pcoord}). Let $K'$ be as in Equation (\ref{eqn:K'}).
    Let $S$ be a sample covariance matrix for $\mathcal{M}_\T$ and let $S'$ be as in  \Cref{eq:S'}.
    This transformation was chosen so that $\mbox{tr}(SK)=\mbox{tr}(S'K')$. Indeed, we can compare the coefficients of each $s_{ij}$ in~$\mbox{tr}(SK)$ and~$\mbox{tr}(S'K')$, denoted $\mathrm{coeff}(s_{ij}, \mbox{tr}(SK))$ and $\mathrm{coeff}(s_{ij}, \mbox{tr}(S'K'))$ respectively. For the diagonal entries corresponding to leaves $i \in \mathrm{Lv}(\T) \setminus \{r\}$, we have
    \[
    \mathrm{coeff}(s_{ii}, \mbox{tr}(S'K')) = \sum_{\substack{t=0 \\ t \neq i}}^n p_{it} = \mathrm{coeff}(s_{ii}, \mbox{tr}(SK)).
    \]
    For $s_{rr}$, we have
    \begin{align*}
        \mathrm{coeff}(s_{rr}, \mbox{tr}(S'K')) &= \sum_{t=1}^n p_{0t} + \sum_{\substack{i=1 \\ i \neq r}}^n \sum_{\substack{t=0\\ t\neq i}}^n p_{it} - 2\sum_{j=1}^n p_{0j} - 2 \sum_{\substack{i < j \\ i, j \neq 0,r}} p_{ij}\\
        &= 2 \sum_{\substack{i = 0 \\i \neq r}}^n \sum_{\substack{j=0 \\j \neq i,r}}^n p_{ij} + \sum_{\substack{j=0\\ j \neq r}}^n p_{rj} - 2 \sum_{\substack{i = 0 \\i \neq r}}^n \sum_{\substack{j=0 \\j \neq i,r}}^n p_{ij} \\
        &= \sum_{\substack{j=0\\ j \neq r}}^n p_{rj} = \mathrm{coeff}(s_{rr}, \mbox{tr}(S'K')).
    \end{align*}
    Similar algebraic manipulations show that the coefficients of $s_{ij}$ for $i \neq j$ are also equal. Hence, the traces are the same.

    By \Cref{thm:detK-formula}, the term $\log \det (K)$ in the log-likelihood function depends only on the unrooted topology of the tree. So,
    \begin{align}\label{eq:mleq S S'}
        \ell_\T(\theta | S)= \ell_{\T'}(\theta | S').
    \end{align} 
    The map sending $S$ to $S'$ is invertible.
    Since the number of complex critical points of the log-likelihood function is fixed and equal to the ML-degree for a generic choice of sample covariance matrix, (a) follows.  

To prove (b), we first claim that the map in \Cref{eqn:K'} sending $K$ to $K'$ and its inverse sending $K'$ to $K$ map the positive definite cone to itself. 
We prove this using Sylvester's Criterion for positive definiteness \cite[Theorem~7.2.5]{horn2012matrix}. 
Without loss of generality, we may let $r = 1$; indeed, if $r \neq 1$, we may conjugate $K'$ by a permutation matrix so that the $r$th row and column become the first. 

For any $S \subset [n]$, let $K_S$ and $K'_S$ be the principal submatrices of $K$ and $K'$ respectively whose rows and columns are indexed by $S$.
For each $i$ such that $0 \leq i \leq n-1$, let $S_i = \{n, n-1, \dots n-i\}$. 
Suppose that $K$ is positive definite. Then by Sylvester's Criterion, $\det(K_{S_i}) > 0$ for each $i$. Since $r = 1$, we have $K_{S_i} = K'_{S_i}$ for each $i \leq n-2$.
So $\det(K'_{S_i}) > 0$ for $i \leq n-2$ as well. Moreover, by the Matrix-Tree Theorem \cite{kirchoff}, we have $\det(K) = \det(K')$ since both $K$ and $K'$ are maximal principal submatrices of the weighted Laplacian of a complete graph on $n+1$ vertices where the weight of the edge $\{i,j\}$ is $p_{ij}$. 
Since $K$ is positive definite, $\det(K) = \det(K') = \det(K'_{S_{n-1}}) > 0$. 
Hence, we have found a nested sequence of length $n$ of principal submatrices of $K'$ with positive determinant. 
By Sylvester's Criterion, $K'$ is positive definite, as needed.

Finally, by \Cref{eq:mleq S S'} and the fact that the map defined by \Cref{eqn:K'} is an automorphism of the positive definite cone, we have that if $K$ maximizes $\ell_\T(\theta \mid S)$ over the positive definite cone, then $K'$ maximizes $\ell_\T(\theta \mid S')$ over the positive definite cone. Hence $K'$ is the MLE for $S'$ in the model $\mathcal{M}_{\T'}$.
\end{proof}}

\section{ML-degrees of BMT models on Star Trees} \label{sec:StarTree}

Let $\T_n$ be the star tree on leaves $\{ 0,1,\ldots,n \}$\textcolor{blue}{,} with unique internal node $x$. We prove \Cref{thm:main-B}, which states that the maximum likelihood degree of its associated BMT model is $2^{n+1}-2n-3$.  
\iftoggle{EA}{}{We will use B\'ezout's Theorem.  
\begin{theorem}[B\'ezout's Theorem {\cite[\S II.2]{basicAGshafarevich}}]
    \label{thm:Bezout}
    Given $n$ hypersurfaces of degrees $d_1, \ldots, d_n$ in a projective space of dimension $n$ over an algebraically closed field, if the intersection of the hypersurfaces is zero-dimensional, then the number of intersection points, counted with multiplicity, is equal to the product of the degrees $d_{1}\cdots d_{n}$.
\end{theorem}
For ease of notation, denote by $\theta_i$ the parameter for edge $\{i,x\}$ in the path parametrization of $\cL_\T^{-1}$. Let $K_{\T_n}(\theta)$ denote the concentration matrix for the BMT model on the star tree $\T_n$ in the path parameters $\theta_0,\dots,\theta_n$. In other words, it is the $n\times n$ matrix whose $(i,j)$ entry is $-\theta_i \theta_j$ when $i \neq j$ and whose $i$th diagonal entry is
\[
\theta_i \sum_{\substack{j=0 \\ j \neq i}}^n \theta_j.
\]  
Let $S\in \PSD^n(\mathbb R)$ be \textcolor{blue}{a} sample covariance matrix. We start by setting up the system of score equations of $\ell_{\T_n}(\theta|S)$. Then we count the solutions $\theta\in \CC^{n+1}$ that have $\det (K_{\T_n}(\theta)) \neq 0$ with their multiplicities. 

\begin{prop} \label{prop:ell T_n}
The score equations of $\ell_{\T_n}(\theta |S)$ have the form 
  \begin{align}\label{eq:score eq}
  \dfrac{\partial \ell_{\T_n}(\theta|S)}{\partial \theta_i} = \frac{1}{\theta_i} + \dfrac{n-1}{\theta_0 + \theta_1 + \cdots + \theta_n} - \sum_{\substack{j = 0 \\ j \neq i}}^n c_{ij} \theta_j, \text{ for } i=0,\ldots, n, 
    \end{align}  where $c_{0j} = s_{jj}, \text{ and } c_{ij} = s_{ii} + s_{jj} - 2 s_{ij} \text{ for } i > 0.$
\end{prop}}

\iftoggle{EA}{\begin{proof}[Proof outline]
    Expand $\mbox{tr} (SK_{\T_n}(\theta))$ and $\det (K_{\T_n}(\theta))$ (using \Cref{thm:detK-formula}).
\end{proof}}%
{\begin{proof}
    Expanding out the expression for the trace of $SK_{\T_n}(\theta)$  and applying \Cref{thm:detK-formula} to~$\T_n$ gives
        \begin{equation*}
      \mbox{tr}(SK_{\T_n}(\theta)) = \sum_{i=0}^n \sum_{j > i} c_{ij} \theta_i \theta_j \text{ and }  \det (K_{\T_n}(\theta))= \prod_{i=0}^n \theta_i \left( \sum_{i=0}^n \theta_i \right)^{n-1}.
        \end{equation*}
    Substituting these expressions into $\ell_{\T_n}(\theta|S)$ and taking its partial derivatives gives exactly the equations in~(\ref{eq:score eq}).
\end{proof}}
For a system $F$  of polynomials, we denote by $V(F)$ its complex affine variety.
In the following steps, we will restate the problem of computing the ML-degree as counting solutions to a polynomial system $F_n$ (\Cref{lm:V(F_n)}) and then as counting solutions to a \emph{homogeneous} polynomial system $\Tilde F_n$ in projective space (\Cref{lm:Tilde F_n}). In all the steps, we need to remove solutions for which $\det (K(\theta))=0$. 

We introduce a new variable $\psi$ which plays the role of the inverse of $\theta_0 + \theta_1 + \cdots + \theta_n$ by adding the equation
 $1- \psi\sum_{i=0}^n \theta_i = 0 $
to the set of the score equations.
For a fixed sample covariance matrix $S=(s_{ij})$ and values $c_{ij}$ as in (\ref{eq:score eq}), let $F_n=\{f_0,\ldots,f_{n+1}\}$ be the following system of $n+2$ polynomials in $\CC[\theta_0,\ldots,\theta_n,\psi]$:
\begin{align}\label{eq:f_s}
    F_n\colon \ {f}_i = 1 + \theta_i \left( (n-1) \psi - \sum_{j \neq i} c_{ij} \theta_j \right), \text{ for } i=0,\ldots,n  \  \text{ and } \ 
{f}_{n+1} = 1 - \psi \left( \sum_{i=0}^n \theta_i \right).
\end{align}
Note that for $i=0,\ldots,n$, the polynomial $f_i$ is  $\dfrac{\partial \ell_{\T_n}(\theta | S)}{\partial \theta_i}$ with its denominator cleared. The goal of Lemma \ref{lm:V(F_n)} is to show that the degree of the affine variety of this system intersected with the algebraic torus is exactly \textcolor{blue}{twice} the maximum likelihood degree of $\cM_{\T_n}$.    
Let $\CC^*$ denote the complex numbers without zero. 

\begin{lemma} \label{lm:V(F_n)} The maximum likelihood degree of the BMT model on $\T_n$ is {half of} the degree of the ideal generated by $F_n$; that is
   { $\mld({\mathcal M}_{\T_n}) = \frac{1}{2} \deg V(F_n)$.}
\end{lemma}

\begin{proof}
    Let $L_n$ denote the ideal generated by the score equations in the ring 
    \[A_n = \CC [\theta_0^\pm, \ldots, \theta_n^\pm, \left( \theta_0 + \dots + \theta_n \right)^{-1}].\] 
    { Each point of $V(L_n)$ corresponds to a critical point of the log-likelihood function. By Theorem \ref{thm:degree}, this correspondence is 2-to-1 since the degree of the parametrization $\phi^\ast_{\T_n}$ is 2. Hence by Corollary \ref{cor:F(k) F(theta)} the  ML-degree of $\cM_{\T_n}$ is $\frac{1}{2} \deg V(L_n)$.}
    We use the variable $\psi$ to represent $\left( \sum_{i = 0}^n \theta_i \right)^{-1}$. By the first isomorphism theorem, 
    \[A_n / L_n \cong \CC [\theta_0^\pm, \ldots, \theta_n^\pm, \psi] / \left( L_n + \langle 1 - \psi \left( \theta_0 + \dots + \theta_n \right) \rangle \right).
    \] By clearing denominators in the score equations, we see that the vanishing locus of $L_n + \langle 1 - \psi \left( \sum_{i=0}^n \theta_i \right) \rangle$ is isomorphic to that of the saturated ideal, $\langle F_n \rangle : \left(\det (K_{\T_n}) \right)^{\infty}$ in $\CC [\theta_0, \ldots, \theta_n, \psi]$. Using the factorization of $\det (K_{\T_n})$ given in \Cref{prop:star-detK-formula}, we see that this is the same as $V(F_n) \cap \left( \CC^\ast \right)^{n+2}$. But if $\theta_i = 0$, then $f_i = 1 \neq 0$, and if $\psi = 0$, then $f_{n+1} = 1 \neq 0$.  It follows that $V(F_n) \cap \left( \CC^\ast \right)^{n+2} = V(F_n)$. Thus, {$\mld({\mathcal M}_{\T_n}) = \frac{1}{2} \deg V(L_n) = \frac{1}{2}\deg V(F_n)$.}
\end{proof}

In order to apply B\'ezout's Theorem, we consider the homogenization of the system~$F_n$. Given a system $G$ of homogeneous polynomials in $m$ variables, let $X(G) \subset \mathbb{P} \CC^{m-1}$ denote its complex projective variety.
The homogenization of the  system  $F_n$ is the system $\Tilde F_n$ of $n+2$ homogeneous polynomials in the $n+3$ variables $\theta_0,\ldots,\theta_n,\psi,z$:    
\begin{align}\label{eq:f_tildas}
   \tilde F_n \colon \     \tilde{f}_i = z^2 + \theta_i \left( (n-1) \psi + \sum_{j \neq i} c_{ij} \theta_j \right), \text{ for } i=0,\ldots,n,   \text{ and } \ \ 
        \tilde{f}_{n+1} = z^2 - \psi \sum_{i=0}^n \theta_i.
\end{align}
Its solution set $X(\Tilde F_n)$ lives in $(n+2)$-dimensional complex projective space. In the next lemma, we consider the {following points in $\mathbb P\CC^{n+2}$, which we call the \emph{standard points}}: \[e_{i}=X(\langle z, \psi, \theta_j \mid j \neq i \rangle) \text{ for } i=0,\ldots,n, \ \text{ and } e_{n+1}=X(\langle z, \theta_i \mid i=0,\dots,n\rangle).\]
We show that the maximum likelihood degree of $\cM_{\T_n}$ is exactly \textcolor{blue}{one-half} the degree of the projective variety $X(\tilde{F}_n)$ with these standard points removed.

\begin{lemma}
    \label{lm:V(Tilde F_n)} 
    For a generic sample covariance matrix, the number of affine solutions to $F_n$, counted with multiplicity, in the algebraic torus is equal to the number of projective solutions, counted with multiplicity, to $\tilde{F}_n$ that are not standard points; that is,
 $\deg \left( V(F_n)\cap {(\CC^{*})}^{n+2} \right) = \deg \left( X(\tilde F_n)\setminus \{e_0,\ldots,e_{n+1}\}\right)$.  
\end{lemma}

\begin{proof}
  We  prove the lemma by showing that    \[X(\Tilde F_n)\setminus (\mathbb P\CC^*)^{n+2}= X(\Tilde F_n,z) =\{e_0,\ldots,e_{n+1}\}.\]   
The bijection 
\((\CC^*)^{n+2} \xrightarrow{ }  (\mathbb P\CC^*)^{n+2} \) sending \((\theta_0,\ldots,\theta_n,\psi)\mapsto{} [\theta_0:\cdots:\theta_n:\psi:1] \)
concludes the rest. 

{
First, we specify the genericity conditions that we impose on the sample covariance matrix. We define the matrix
\[
C = \begin{pmatrix}
        0 & c_{01} & \cdots & c_{0n} & n-1 \\
        c_{01} & 0 & \cdots & c_{1n} & n-1 \\
        \vdots & \vdots & \ddots & \vdots & \vdots \\
        c_{0n} & c_{1n} & \cdots & 0 & n-1 \\
        n-1 & n-1 & \cdots & n-1 & 0
    \end{pmatrix}.
\]
The entries of $C$ are an invertible linear function of the entries of the sample covariance matrix. So we may specify these genericity conditions on $C$. Let the entries of $C$ be such that  the determinant of every principal submatrix of $C$ of size greater than 1 is non-zero. Each of these is a non-trivial polynomial condition on the entries of $C$, so a generic $C$ satisfies these constraints.

Take $P=[\theta_0:\cdots:\theta_n:\psi:z]\in X(\Tilde F_n)\setminus (\mathbb P\CC^*)^{n+2}$. Points not in the torus $(\PP\CC^\ast)^{n+2}$ have at least one coordinate equal to zero. If any of the $\theta_i$ are zero, then $\tilde{f}_i (P) = z^2$, so $z = 0$ as well. When $\psi = 0$, we have $\tilde{f}_{n+1}(P) = 0$, so again $z = 0$. So in all cases, the $z$ entry of $P$ is zero.  Hence the system $\tilde{F}_n$ is  \begin{gather*}
        \tilde{f}_i(P) = \theta_i \left( (n-1) \psi + \sum_{j \neq i} c_{ij} \theta_j \right) \text{ for } i=0,\ldots,n, \text{ and }
        \tilde{f}_{n+1}(P) = \psi \left( \sum_{i=0}^n \theta_i \right).
    \end{gather*}
    
Let $A \subset \{0,\dots,n+1\}$ be the set of indices of non-zero coordinates of $P$, where the index $n+1$ is associated with $\psi$. For the sake of contradiction, suppose that $\#A > 1$. If $i \in A$, then in order for $\tilde{f}_i(P)$ to vanish, we must have
\begin{gather*}
    (n-1)\psi + \sum_{\substack{j \neq i \\ j \in A}} \theta_j = 0 \ \text{ if } \ i =0,\ldots,n, \ \text{ and } \ \sum_{\substack{j \leq n \\ j \in A}} \theta_j {\color{blue}=0}\ \text{ if } \ i = n+1.
\end{gather*}
 If such a solution exists, then the $A \times A$ minor of $C$ has a non-trivial kernel\footnote{\color{blue}Note that if $n+1 \not\in A$, then the $\psi$ term of these equations is zero.}. But this contradicts our assumption that all principal submatrices of $C$ of size greater than 1 are invertible. So we must have $\#A = 1$ and $P$ must be equal to one of $e_0,\dots,e_{n+1}$. Moreover, each of these clearly lies in $X(\tilde{F}_n)$, as needed.}
\end{proof}

To finish our proof, we must compute the multiplicities of points $e_0,\ldots,e_{n+1}$ in $X(\Tilde F_n)$, and subtract them from $\deg X(\Tilde F_n)$. Let $\CC[\![ x_0,\ldots,x_n]\!]$ denote the ring of formal power series in variables~$x_0,\dots,x_n$ {, which is the completion of $\CC[ x_0,\ldots,x_n]$ at the origin}. In order to compute these multiplicities, we make use of \Cref{thm:dim}. The \emph{standard monomials} of an ideal $I$ with respect to a local order are the monomials $x^\alpha$ such that $x^\alpha$ does not belong to the leading term ideal $\mathrm{LT}(I)$ with respect to this order.

\begin{theorem}[\cite{usingAG}, \S 4.4, Theorem~4.3]
\label{thm:dim} Let $\widehat R = \CC[\![x_0,\ldots,x_n]\!]$. 
Let  $\widehat J\subset \widehat R$ be an ideal, $>$ a local order, and  $\mathrm{LT}(\widehat J)$ the leading term ideal for $\widehat J$ with respect to $>$. 
If $\widehat R/\widehat J$ contains finitely many standard monomials,  
then  $\dim_\CC(\widehat R/\widehat J)$ is the number of standard monomials.
\end{theorem}

{
The following lemma draws a connection between the multiplicity of the origin in a variety and the dimension of the quotient in the ring of formal power series.

\begin{lemma}
    \label{lm:helper}
    Let $I$ be an ideal in $\CC[x_1,\ldots,x_n]$, and let $p = \langle x_1, \ldots, x_n \rangle$. Then
    \begin{equation}
        \label{eqn:helper}
        \dim_\CC (\CC[x_1,\ldots,x_n]/I)_p = \dim_\CC \CC[\![ x_1, \ldots, x_n ]\!]/I \CC[\![ x_1, \ldots, x_n ]\!],
    \end{equation}
    provided that at least one of the above quantities is finite.
\end{lemma}

\begin{proof}
    Let $A = (\CC[x_1,\ldots,x_n]/I)_p$. Note that $A$ is a Noetherian local ring.
    Then \Cref{eqn:helper} reads $\dim_\CC A = \dim_\CC \widehat A$, where $\widehat A$ is the completion of $A$.
    Our first step is to show that if either of the quantities in \Cref{eqn:helper} is finite, then the Krull dimensions $\dim A = \dim \widehat A = 0$. 
    We will use the following two fundamental facts from commutative algebra:
    \begin{enumerate}
        \item In a Noetherian local ring, completion preserves Krull dimension, so $\dim A = \dim \widehat A$.
        \item If a $\CC$-algebra $R$ is finite-dimensional as a $\CC$-vector space, then $\dim R = 0$.
    \end{enumerate}
    Combining the two facts above, if $\dim_\CC A$ is finite, then $\dim A = 0$, and so $\dim \widehat A = 0$.
    On the other hand, if $\dim_\CC \widehat A$ is finite, then $\dim \widehat A = 0$, and so $\dim A = 0$.
    
    Let $I = \bigcap_{i=1}^k q_i$ be a primary decomposition of $I$, and let $J$ be the intersection over the primary components $q_i$ contained in the maximal ideal $p$.
    The dimension of $A$ is the maximum of the dimensions of all primary components of $I$ contained in $p$; hence, $J$ is zero-dimensional.
    Moreover, we have $A = (\CC[x_1,\ldots,x_n]/J)_p$.
    Now apply Proposition 2.11 of \cite{usingAG} to $J$, which proves \Cref{eqn:helper} for zero-dimensional ideals.
\end{proof}}

In the proof of the next lemma, we define a local order on a given power series ring. We set a variable $\theta_i$ or $\psi$ equal to 1 to localize at the prime ideal of the corresponding standard point $e_i$. Then we use the polynomials $\tilde f_0, \dots, \tilde f_{n+1}$ to expand each $\theta_i$ and $\psi$ as a power series in $z$. This allows us to find the standard monomials and compute the multiplicity $e_i.$

\begin{lemma}
\label{lm:Tilde F_n}
\phantom{hello world} 

\begin{enumerate}
    \item The multiplicity of the standard point $e_i$ for $i=0,\ldots, n $ in $X(\tilde{F}_n)$ is~four.
    \item The multiplicity of the standard point $e_{n+1}$ in $X(\tilde{F}_n)$ is two.
    \end{enumerate}
\end{lemma}

\begin{proof}[Proof of (1)]
 By symmetry, we only need to prove the lemma for $e_0$. Let $R=\CC[\theta_1, \ldots, \theta_n, \psi, z] $ and  $\widehat{R}=\CC[\![\theta_1, \ldots, \theta_n, \psi, z]\!]$.  Substituting $\theta_0=1$ to \Cref{eq:f_tildas}, we obtain
    \begin{align*}
        &  \bar{f}_0 = z^2 +   (n-1)\psi  + \sum_{j \neq 0} c_{0j} \theta_j, \\
        &    \bar{f}_i = z^2 + \theta_i \left( (n-1)\psi  + c_{0i}+ \sum\limits_{j \neq i,0} c_{ij} \theta_j \right) \ \text{ for } i=1,\ldots,n, \text{ and }\\
        &  \bar{f}_{n+1} = z^2 - \psi\left( 1+\sum\limits_{i=1}^n \theta_i \right).
    \end{align*}
    Denote $J = \langle \bar{f}_0, \ldots, \bar{f}_{n+1} \rangle\subseteq R$ and $p=\langle \theta_1,\ldots, \theta_n,z,\psi \rangle\subseteq R$.  
    Let $\left(R / J \right)_p$ be the local ring at $p$. 
    By definition, the intersection multiplicity of $e_0$ is { $\dim_\CC(R / J)_p$. By Lemma \ref{lm:helper}, this is equal to $\dim_\CC(\widehat R/\widehat J)$ as long as either one of these dimensions is finite.
    Here, $\widehat J$ is the completion of $J$ at the origin.}
    By \Cref{thm:dim}, the $\dim_{\CC} (\widehat R /\widehat J)$ is equal to the number of standard monomials of $\widehat J$ with respect to any local order, provided that either number is finite.
    In order to facilitate the computation, we find a standard basis for $\widehat J$.

    First, we use the functions $\bar{f}_1,\ldots,\bar{f}_{n+1}$  to write power series expansions for $\psi$ and $\theta_j$ for  $j = 1, \ldots, n$  in terms of $z$. Since  $\bar{f}_j = 0$ for $j=1,\dots,n$, we have
    \[
    \theta_j = \frac{-\frac{1}{c_{0j}} z^2}{1 + (n-1) \psi + \sum\limits_{i \neq 0, j} \frac{c_{ij}}{c_{0j}} \theta_i} = -\frac{1}{c_{0j}} z^2 \left(1 - \frac{n-1}{c_{0j}} \psi - \sum_{i \neq 0,j} \frac{c_{ij}}{c_{0j}} \theta_i + \cdots \right).
    \]
    Similarly, since $\bar{f}_{n+1} = 0$, we have
    \[
    \psi = \frac{z^2}{1 + \sum\limits_{i \neq 0} \theta_i} = z^2 \left( 1 - \sum\limits_{i \neq 0} \theta_i + \cdots \right).
    \]
    By substitution, we obtain power series expansions for $\theta_j$, $j \neq 0$, and $\psi$ up to degree 4 in $z$.
    \begin{align}\label{eq:zetas}
        \theta_j = -\frac{1}{c_{0j}} z^2 + \left( - \frac{n-1}{c_{0j}} + \sum_{i \neq 0,j} \frac{c_{ij}}{c_{0j}^2} \right) z^4 + O(z^6) \text{ and }
        \psi = z^2 + \left( \sum_{i \neq 0} \frac{1}{c_{0i}} \right) z^4 + O(z^6).
    \end{align}
    Denote by $\bar{g}_0$ the equation obtained by writing  all the variables  in $\bar f_0$ in terms of $z$ using the equations in (\ref{eq:zetas}), so that  {\small 
    \begin{align*}
        \bar{g}_0 &= z^2 + (n-1) \left( z^2 + \left( \sum_{i \neq 0} \frac{1}{c_{0i}} \right) z^4 + O(z^6) \right) \\
        & \phantom{=} + \sum_{i \neq 0} c_{0i} \left( -\frac{1}{c_{0i}} z^2 + \left( - \frac{n-1}{c_{0i}} + \sum_{k \neq 0,i} \frac{c_{ik}}{c_{0i}^2} \right) z^4 + O(z^6) \right) \\
        &= \left( \sum_{i \neq 0} \frac{1}{c_{0i}} \left( 1 + \sum_{k \neq 0, i} c_{ik} \right) - n(n-1) \right) z^4 + O(z^6).
    \end{align*}}
    
By the genericity of the $c_{ij}$'s, the coefficient of $z^4$ in the power series $g_0$ is generically non-zero. Let~$G=\{\bar{g}_0, \bar{f}_1, \ldots, \bar{f}_n, \bar{f}_{n+1}\}$. Note that $\widehat J = \langle G \rangle$. We further claim that $G$ is a standard basis for~$\widehat J$. 
    Indeed, take  $<$ to be a negative graded monomial order, i.e.
    \(1 > \theta_1, \ldots, \theta_j, \psi, z > \theta_i \theta_j, \ldots\), and so on. 
    Since the $<$-leading terms of $\bar{g}_0$, $\bar{f}_j, \ldots, \bar{f}_{n+1}$ are relatively prime, the set $G$ is a standard basis for $\widehat{J}$. Thus, $\mbox{LT}(\widehat J) = \langle \theta_1, \ldots, \theta_n, \psi, z^4 \rangle$. There are four standard monomials, $1, z, z^2, z^3$, and thus $\dim_\CC \widehat R/ \widehat J = 4$ by \Cref{thm:dim}. 
    Since $\dim_\CC \widehat R/ \widehat J = \mathrm{length}(R/J)_p = \mathrm{mult}(e_0)$, this implies $\mathrm{mult}(e_0)=4$.
\end{proof}

\begin{proof}[Proof of (2)]
    After dehomogenizing by $\psi = 1$, we obtain the system
    \begin{align*}
        &  \bar{f}_i = z^2 + \theta_i \left( (n-1)  + \sum\limits_{j \neq i} c_{ij} \theta_j \right) \ \text{ for } i=0,1,\ldots,n,\\
        &  \bar{f}_{n+1} = z^2 - \left( \sum\limits_{i=0}^n \theta_i \right).
    \end{align*}
    As in the previous case, let $R= \CC[\theta_0, \ldots, \theta_n, z]$,  $J = \langle \bar{f}_0, \ldots, \bar{f}_n, \bar{f}_{n+1} \rangle$ and $p=\langle \theta_0,\ldots, \theta_n,z\rangle $. 
    By \Cref{lm:helper}, it suffices to compute $\mathrm{dim}_\CC(\widehat{R}/\widehat{J})$.
    We begin by finding a standard basis for $\widehat J$.  
    Since $\bar{f}_i = 0$ for $i=0,\dots,n$, we may solve for $\theta_i$ in terms of $z$ and obtain
    \begin{align*}
        \theta_i = \frac{-z^2}{(n-1) + \sum\limits_{j \neq i} c_{ij} \theta_j}
        = - \frac{z^2}{n - 1} \left( 1 - \sum\limits_{j \neq i} c_{ij} \theta_j \right) 
        = z^2 \left( - \frac{1}{n - 1} + O(z^4) \right).
    \end{align*}
    Substituting the relations derived above into $\bar{f}_{n+1}$, the power series, 
    $$\bar{g}_{n+1} = z^2 + \frac{n + 1}{n - 1} z^2 + O(z^4) = \frac{2n}{n - 1} z^2 + O(z^4).$$
    Note that $\widehat{J} = \langle \bar{f}_0, \ldots, \bar{f}_n, \bar{g}_{n+1} \rangle$. For any negative graded monomial order $<$, the leading terms of $\bar{f}_0, \ldots, \bar{f}_n, \bar{g}_{n+1}$ are relatively prime. Thus, $\bar{f}_0, \ldots, \bar{f}_n, \bar{g}_{n+1}$ form a standard basis of $\widehat{J}$. It follows that 
    $$ \mbox{LT}(\widehat J) = \langle \theta_0, \ldots, \theta_n, z^2 \rangle.$$
  By \Cref{thm:dim}, $\dim_\CC \widehat R/ \widehat J = 2$, and thus $\mathrm{mult}(e_{n+1})=2$. 
\end{proof}
{
\begin{cor}\label{cor:zerodim}
    The variety $X(\tilde{F}_n)$ is zero-dimensional.
\end{cor}}
{
\begin{proof}
    For the sake of contradiction, suppose that $X(\tilde{F}_n)$ contains a positive-dimensional irreducible component $W$, and let $H = V(z)$ be the hyperplane at infinity.
    Since $H$ is a hypersurface and the dimension of $W$ is at least one, their intersections are nonempty. Let $p$ be a point in $H \cap W$.
    By \Cref{lm:V(Tilde F_n)}, $p \in \{ e_0, \ldots, e_{n+1} \}$.
    
    On the other hand, \Cref{lm:Tilde F_n} shows that for any $p \in \{ e_0, \ldots, e_{n+1} \}$, $(R/I)_p$ is a finite-dimensional $\CC$-vector space.
    Thus, the Krull dimension of $(R/I)_p$ is zero. 
    Since this is the maximum of the dimensions of all irreducible components of $V(I)$ containing $p$, $W$ must have dimension zero.
    This contradicts the assumption that $W$ is positive-dimensional.
    Therefore, $X(\tilde{F}_n)$ contains no positive-dimensional components and hence is zero-dimensional.
\end{proof}}

We are ready to prove the main result.

\begin{reptheorem}{thm:main-B}
    \label{thm:star_tree_MLD}
    The maximum likelihood degree of the Brownian motion star tree model on $n+1$ leaves is $2^{n+1}-2n-3$. 
\end{reptheorem}
 
\iftoggle{EA}{\begin{proof}[Proof outline]
    We begin by adding a new variable to our system of score equations, $\psi := \left( \sum_{i = 0}^n \theta_i \right)^{-1}$, along with a new constraint, $f_\psi := 1 - \psi \left( \theta_0 + \theta_1 + \cdots + \theta_n \right)$. After clearing denominators and homogenizing, this becomes the system $\tilde{F}_n$ consisting of equations
    \begin{align}\label{eq:f_tildas}
             \tilde{f}_i := z^2 + \theta_i \left( (n-1) \psi + \sum_{j \neq i} c_{ij} \theta_j \right), \text{ for } i=0,\ldots,n,   \text{ and } \ \ 
            \tilde{f}_{n+1} := z^2 - \psi \sum_{i=0}^n \theta_i.
    \end{align}
    We show that the critical points of $\ell_{\T_n}(K_{\T_n}(\theta))$ in $\mathcal{L}_{\T_n}^{-1}$ are in bijection with $V(\tilde{F}_n) \cap \left( \CC^\ast \right)^{n+3}$. The points of $V(\tilde{F}_n) \cap \left( \CC^\ast \right)^{n+3}$ can be counted by first using B\'ezout's Theorem to show that $\deg V(\tilde{F}_n) = 2^{n+2}$ and then counting the non-torus points with multiplicity. We find that for generic $S$, the only non-torus points are the standard basis vectors  $e_0,\cdots,e_{n+1}$ in $\PP\CC^{n+2}$ (but not $e_{n+2}$). Using local orders, we prove that the multiplicity of $e_{n+1}$ in $V(\tilde{F}_n)$ is $2$ and that the multiplicity of $e_i$  is $4$ for $i = 0, \ldots, n$. Thus, $\ell_{\T_n}(K_{\T_n}(\theta))$ has $2^{n+2} - 4(n+1) - 2$ critical points. By \Cref{cor:F(k) F(theta)}, $\ell_{\T_n}(K_{\T_n})$ has  $\frac12 \left( 2^{n+2} - 4(n+1) - 2 \right)$ critical points, counting multiplicity. 
\end{proof}}%
{\begin{proof}
    { The system $\Tilde F_n$ is a homogeneous system of $n+2$ quadratic polynomials whose solution set is zero-dimensional by \Cref{cor:zerodim}. So by B\'ezout's Theorem, it has $2^{n+2}$ solutions, counted with multiplicity.} By \Cref{lm:V(F_n)} and \Cref{lm:V(Tilde F_n)}, we have
    \[
     \mld(\mathcal M_{\T_n})=\dfrac{1}{2}\deg \left( X(\tilde F_n)\setminus \{e_1,\ldots,e_{n+2}\}\right).
    \]
    Applying \Cref{lm:Tilde F_n} to remove these standard points with their multiplicities, we obtain
    \begin{align*}
        \mld(\mathcal M_{\T_n})&=\dfrac{1}{2}\left(\deg  X(\tilde F_n) -\sum_{i=1}^{n+2}\mathrm{mult}(e_i)\right)\\
        &=\dfrac{1}{2}\left(2^{n+2}-4(n+1)-2)\right)\\
        & =2^{n+1}-2n-3. \qedhere
    \end{align*}
\end{proof}}

\section{Discussion}
\label{sec:Discussion}
In this paper, we use algebraic techniques to give a formula for the ML-degree of the BMT model on a star tree. Theorem \ref{thm:detK-formula} is a generalization of the Cayley-Pr\"ufer Theorem, which gives a formula for the determinant of a matrix in $\mathcal{L}^{-1}_\T$. We used this result to show that the ML-degree of the BMT model is the same for all trees with the same unrooted topology. 


\begin{figure}[!htpb]{
        \centering
        \begin{tabular}{c|ccc||c|ccc}
 tree topology  & 
deg &   rmldeg &   mldeg &
 tree topology  & 
deg &   rmldeg &   mldeg
\\   
\hline
\hline
 
 \begin{tabular}{c}\scalebox{0.4}{\begin{tikzpicture}
  [scale=.6,auto=left,very thick,-]
  \node [] (nz) at (3,3.6)  {};
\node [circle, fill=red!60!blue!50!white] (n0) at (3,3)  {};
\node [circle, fill=red!60!blue!50!white] (n1) at (-3,3)  {};
\node [circle, fill=red!60!blue!50!white] (n2) at (-4,0) {};
\node [circle, fill=red!60!blue!50!white](n3) at (-3,-3)  {};
\node [circle, fill=red!60!blue!50!white] (n4) at (4,1) {};
\node [circle, fill=red!60!blue!50!white](n5) at (4,-1)  {};
\node [circle, fill=red!60!blue!50!white](n8) at (3,-3)  {};
\node [circle, fill=gray!50] (n7) at (1.5,0)  {};
\node [circle, fill=gray!50] (n9) at (-1.5,0)  {};
\draw (n0) -- (n7) [green!40!black]  node [midway, right, black] {};
\draw (n9) -- (n1) [green!40!black]  node [midway,left, black  ] {};
\draw (n9) -- (n2) [green!40!black] node [midway, right, black ] {};
\draw (n9) -- (n3) [green!40!black] node [midway,left , black ] {};
\draw (n7) -- (n4) [green!40!black] node [near end,right , black ] {};
\draw (n7) -- (n8) [green!40!black] node [near end,right , black ] {};
\draw (n7) -- (n9) [green!40!black] node [near end,right , black ] {};
\draw (n7) -- (n5) [green!40!black]  node [midway, black] {};
\end{tikzpicture}} \end{tabular}
& $93$ & $44$ &  $259$ &
  \begin{tabular}{c}\scalebox{0.4}{\begin{tikzpicture}
  [scale=.6,auto=left,very thick,-]
  \node [] (nz) at (3,3.6)  {};
\node [circle, fill=red!60!blue!50!white] (n0) at (3,3)  {};
\node [circle, fill=red!60!blue!50!white] (n1) at (-4,2.5)  {};
\node [circle, fill=red!60!blue!50!white] (n2) at (-4,-2.5) {};
\node [circle, fill=red!60!blue!50!white](n3) at (4,-2)  {};
\node [circle, fill=red!60!blue!50!white] (n4) at (3,-3) {};
\node [circle, fill=red!60!blue!50!white](n5) at (5,0)  {};
\node [circle, fill=red!60!blue!50!white](n8) at (4,2)  {};
\node [circle, fill=gray!50] (n7) at (2,0)  {};
\node [circle, fill=gray!50] (n9) at (-2,0)  {};
\draw (n0) -- (n7) [green!40!black]  node [midway, right, black] {};
\draw (n9) -- (n1) [green!40!black]  node [midway,left, black  ] {};
\draw (n9) -- (n2) [green!40!black] node [midway, right, black ] {};
\draw (n7) -- (n3) [green!40!black] node [midway,left , black ] {};
\draw (n7) -- (n4) [green!40!black] node [near end,right , black ] {};
\draw (n7) -- (n8) [green!40!black] node [near end,right , black ] {};
\draw (n7) -- (n9) [green!40!black] node [near end,right , black ] {};
\draw (n7) -- (n5) [green!40!black]  node [midway, black] {};
\end{tikzpicture}} \end{tabular}
& $95$ & $26$ &  $53$ \\
 \hline
  
 \begin{tabular}{c}\scalebox{0.4}{\begin{tikzpicture}
  [scale=.6,auto=left,very thick,-]
  \node [] (nz) at (3,3.6)  {};
\node [circle, fill=red!60!blue!50!white] (n0) at (0,3)  {};
\node [circle, fill=red!60!blue!50!white] (n1) at (-4,3)  {};
\node [circle, fill=red!60!blue!50!white] (n2) at (-5,0) {};
\node [circle, fill=red!60!blue!50!white](n3) at (-4,-3)  {};
\node [circle, fill=red!60!blue!50!white] (n4) at (5,0) {};
\node [circle, fill=red!60!blue!50!white](n5) at (4,-3)  {};
\node [circle, fill=red!60!blue!50!white](n8) at (4,3)  {};
\node [circle, fill=gray!50] (n7) at (0,0)  {};
\node [circle, fill=gray!50] (n9) at (-2.5,0)  {};
\node [circle, fill=gray!50] (n10) at (2.5,0)  {};
\draw (n0) -- (n7) [green!40!black]  node [midway, right, black] {};
\draw (n9) -- (n1) [green!40!black]  node [midway,left, black  ] {};
\draw (n9) -- (n2) [green!40!black] node [midway, right, black ] {};
\draw (n9) -- (n3) [green!40!black] node [midway,left , black ] {};
\draw (n10) -- (n4) [green!40!black] node [near end,right , black ] {};
\draw (n10) -- (n8) [green!40!black] node [near end,right , black ] {};
\draw (n7) -- (n9) [green!40!black] node [near end,right , black ] {};
\draw (n7) -- (n10) [green!40!black] node [near end,right , black ] {};
\draw (n10) -- (n5) [green!40!black]  node [midway, black] {};
\end{tikzpicture}} \end{tabular}
& $90$ & $16$ &  $221$ &
  \begin{tabular}{c}\scalebox{0.4}{\begin{tikzpicture}
  [scale=.6,auto=left,very thick,-]
\node [circle, fill=red!60!blue!50!white] (n0) at (3,3)  {};
\node [circle, fill=red!60!blue!50!white] (n1) at (-8,3)  {};
\node [circle, fill=red!60!blue!50!white] (n2) at (-8,-3) {};
\node [circle, fill=red!60!blue!50!white](n3) at (-4,-3)  {};
\node [circle, fill=red!60!blue!50!white] (n4) at (-4,3) {};
\node [circle, fill=red!60!blue!50!white](n5) at (-2,-3)  {};
\node [circle, fill=red!60!blue!50!white](n6) at (3,-3)  {};
\node [circle, fill=gray!50] (n7) at (0,0)  {};
\node [circle, fill=gray!50] (n8) at (-2,0)  {};
\node [circle, fill=gray!50] (n9) at (-4,0)  {};
\node [circle, fill=gray!50] (n10) at (-6,-0)  {};
\draw (n0) -- (n7) [green!40!black]  node [midway, right, black] {};
\draw (n7) -- (n6) [green!40!black]  node [midway,left, black  ] {};
\draw (n7) -- (n8) [green!40!black] node [midway, right, black ] {};
\draw (n8) -- (n5) [green!40!black] node [midway,left , black ] {};
\draw (n8) -- (n9) [green!40!black] node [near end,right , black ] {};
\draw (n9) -- (n4) [green!40!black] node [near end,right , black ] {};
\draw (n9) -- (n3) [green!40!black] node [near end,right , black ] {};
\draw (n9) -- (n10) [green!40!black]  node [midway, black] {};
\draw (n2) -- (n10) [green!40!black]  node [midway, black] {};
\draw (n1) -- (n10) [green!40!black]  node [midway, black] {};
\end{tikzpicture}} \end{tabular}
& $51$ & $4$ &  $83$ \\
 \hline

 \begin{tabular}{c}\scalebox{0.4}{\begin{tikzpicture}
  [scale=.6,auto=left,very thick,-]
  \node [] (nz) at (3,3.8)  {};
  \node [circle, fill=red!60!blue!50!white] (n0) at (4,3)  {};
\node [circle, fill=red!60!blue!50!white] (n1) at (-4,3)  {};
\node [circle, fill=red!60!blue!50!white] (n2) at (-5,0) {};
\node [circle, fill=red!60!blue!50!white](n3) at (-4,-3)  {};
\node [circle, fill=red!60!blue!50!white] (n4) at (0,3) {};
\node [circle, fill=red!60!blue!50!white](n5) at (0,-3)  {};
\node [circle, fill=red!60!blue!50!white](n6) at (4,-3)  {};
\node [circle, fill=gray!50] (n7) at (2,0)  {};
\node [circle, fill=gray!50] (n8) at (0,0)  {};
\node [circle, fill=gray!50] (n9) at (-2,0)  {};
\draw (n0) -- (n7) [green!40!black]  node [midway, right, black] {};
\draw (n7) -- (n6) [green!40!black]  node [midway,left, black  ] {};
\draw (n7) -- (n8) [green!40!black] node [midway, right, black ] {};
\draw (n8) -- (n5) [green!40!black] node [midway,left , black ] {};
\draw (n8) -- (n4) [green!40!black] node [near end,right , black ] {};
\draw (n8) -- (n9) [green!40!black] node [near end,right , black ] {};
\draw (n9) -- (n3) [green!40!black] node [near end,right , black ] {};
\draw (n9) -- (n2) [green!40!black]  node [midway, black] {};
\draw (n9) -- (n1) [green!40!black]  node [midway, black] {};
\end{tikzpicture}} \end{tabular}
& $77$ & $16$ &  $181$ &
  \begin{tabular}{c}\scalebox{0.4}{\begin{tikzpicture}
  [scale=.6,auto=left,very thick,-]
\node [circle, fill=red!60!blue!50!white] (n0) at (4,3)  {};
\node [circle, fill=red!60!blue!50!white] (n1) at (-4,3)  {};
\node [circle, fill=red!60!blue!50!white] (n2) at (-4,-3) {};
\node [circle, fill=red!60!blue!50!white](n3) at (0,3)  {};
\node [circle, fill=red!60!blue!50!white] (n4) at (-1,-3) {};
\node [circle, fill=red!60!blue!50!white](n5) at (1,-3)  {};
\node [circle, fill=red!60!blue!50!white](n6) at (4,-3)  {};
\node [circle, fill=gray!50] (n7) at (2,0)  {};
\node [circle, fill=gray!50] (n8) at (0,0)  {};
\node [circle, fill=gray!50] (n9) at (-2,0)  {};
\draw (n0) -- (n7) [green!40!black]  node [midway, right, black] {};
\draw (n7) -- (n6) [green!40!black]  node [midway,left, black  ] {};
\draw (n7) -- (n9) [green!40!black] node [midway, right, black ] {};
\draw (n8) -- (n5) [green!40!black] node [midway,left , black ] {};
\draw (n8) -- (n4) [green!40!black] node [near end,right , black ] {};
\draw (n8) -- (n9) [green!40!black] node [near end,right , black ] {};
\draw (n8) -- (n3) [green!40!black] node [near end,right , black ] {};
\draw (n9) -- (n2) [green!40!black]  node [midway, black] {};
\draw (n9) -- (n1) [green!40!black]  node [midway, black] {};
\end{tikzpicture}} \end{tabular}
& $47$ & $11$ &  $81$ \\
 \hline

 \begin{tabular}{c}\scalebox{0.4}{\begin{tikzpicture}
  [scale=.6,auto=left,very thick,-]
\node [circle, fill=red!60!blue!50!white] (n0) at (5,3)  {};
\node [circle, fill=red!60!blue!50!white] (n1) at (-5,3)  {};
\node [circle, fill=red!60!blue!50!white] (n2) at (-5,0) {};
\node [circle, fill=red!60!blue!50!white](n3) at (-5,-3)  {};
\node [circle, fill=red!60!blue!50!white] (n4) at (-1,-3) {};
\node [circle, fill=red!60!blue!50!white](n5) at (1,-3)  {};
\node [circle, fill=red!60!blue!50!white](n6) at (5,-3)  {};
\node [circle, fill=gray!50] (n7) at (3,0)  {};
\node [circle, fill=gray!50] (n8) at (1,0)  {};
\node [circle, fill=gray!50] (n9) at (-1,0)  {};
\node [circle, fill=gray!50] (n10) at (-3,0)  {};
\draw (n0) -- (n7) [green!40!black]  node [midway, right, black] {};
\draw (n7) -- (n6) [green!40!black]  node [midway,left, black  ] {};
\draw (n7) -- (n8) [green!40!black] node [midway, right, black ] {};
\draw (n8) -- (n5) [green!40!black] node [midway,left , black ] {};
\draw (n8) -- (n9) [green!40!black] node [near end,right , black ] {};
\draw (n9) -- (n4) [green!40!black] node [near end,right , black ] {};
\draw (n10) -- (n3) [green!40!black] node [near end,right , black ] {};
\draw (n9) -- (n10) [green!40!black]  node [midway, black] {};
\draw (n2) -- (n10) [green!40!black]  node [midway, black] {};
\draw (n1) -- (n10) [green!40!black]  node [midway, black] {};
\end{tikzpicture}} \end{tabular}
& $61$ & $4$ &  $115$ &
  \begin{tabular}{c}\scalebox{0.4}{\begin{tikzpicture}
  [scale=.6,auto=left,very thick,-]
  \node [] (nz) at (3,5.6)  {};
\node [circle, fill=red!60!blue!50!white] (n0) at (0,5)  {};
\node [circle, fill=red!60!blue!50!white] (n1) at (-5,4)  {};
\node [circle, fill=red!60!blue!50!white] (n2) at (-5,-3) {};
\node [circle, fill=red!60!blue!50!white](n3) at (-1,-4)  {};
\node [circle, fill=red!60!blue!50!white] (n4) at (1,-4) {};
\node [circle, fill=red!60!blue!50!white](n5) at (5,-3)  {};
\node [circle, fill=red!60!blue!50!white](n6) at (5,4)  {};
\node [circle, fill=gray!50] (n7) at (0,1.8)  {};
\node [circle, fill=gray!50] (n8) at (-2.5,1)  {};
\node [circle, fill=gray!50] (n9) at (0,-.4)  {};
\node [circle, fill=gray!50] (n10) at (2.5,1)  {};
\draw (n0) -- (n7) [green!40!black]  node [midway, right, black] {};
\draw (n7) -- (n10) [green!40!black]  node [midway,left, black  ] {};
\draw (n7) -- (n9) [green!40!black] node [midway, right, black ] {};
\draw (n7) -- (n8) [green!40!black] node [midway,left , black ] {};
\draw (n8) -- (n1) [green!40!black] node [near end,right , black ] {};
\draw (n8) -- (n2) [green!40!black] node [near end,right , black ] {};
\draw (n9) -- (n3) [green!40!black] node [near end,right , black ] {};
\draw (n9) -- (n4) [green!40!black]  node [midway, black] {};
\draw (n5) -- (n10) [green!40!black]  node [midway, black] {};
\draw (n6) -- (n10) [green!40!black]  node [midway, black] {};
\end{tikzpicture}} \end{tabular}
& $42$ & $4$ &  $63$ \\
 \hline

 \begin{tabular}{c}\scalebox{0.4}{\begin{tikzpicture}
  [scale=.6,auto=left,very thick,-]
  \node [circle, fill=red!60!blue!50!white] (n0) at (5,3)  {};
\node [circle, fill=red!60!blue!50!white] (n1) at (-5,3)  {};
\node [circle, fill=red!60!blue!50!white] (n2) at (-5,1) {};
\node [circle, fill=red!60!blue!50!white](n3) at (-5,-1)  {};
\node [circle, fill=red!60!blue!50!white] (n4) at (-5,-3) {};
\node [circle, fill=red!60!blue!50!white](n5) at (0,-3)  {};
\node [circle, fill=red!60!blue!50!white](n6) at (5,-3)  {};
\node [circle, fill=gray!50] (n7) at (2,0)  {};
\node [circle, fill=gray!50] (n8) at (0,0)  {};
\node [circle, fill=gray!50] (n9) at (-2,0)  {};
\draw (n0) -- (n7) [green!40!black]  node [midway, right, black] {};
\draw (n7) -- (n6) [green!40!black]  node [midway,left, black  ] {};
\draw (n7) -- (n8) [green!40!black] node [midway, right, black ] {};
\draw (n8) -- (n5) [green!40!black] node [midway,left , black ] {};
\draw (n8) -- (n9) [green!40!black] node [near end,right , black ] {};
\draw (n9) -- (n4) [green!40!black] node [near end,right , black ] {};
\draw (n9) -- (n3) [green!40!black] node [near end,right , black ] {};
\draw (n9) -- (n2) [green!40!black]  node [midway, black] {};
\draw (n9) -- (n1) [green!40!black]  node [midway, black] {};
\end{tikzpicture}} \end{tabular}
& $60$ & $11$ &  $101$ &
  \begin{tabular}{c}\scalebox{0.4}{\begin{tikzpicture}
  [scale=.6,auto=left,very thick,-]
  \node [] (nz) at (3,3.6)  {};
\node [circle, fill=red!60!blue!50!white] (n0) at (0,3)  {};
\node [circle, fill=red!60!blue!50!white] (n1) at (-4,3)  {};
\node [circle, fill=red!60!blue!50!white] (n2) at (-4,-3) {};
\node [circle, fill=red!60!blue!50!white](n3) at (-1.5,-4)  {};
\node [circle, fill=red!60!blue!50!white] (n4) at (1.5,-4) {};
\node [circle, fill=red!60!blue!50!white](n5) at (4,3)  {};
\node [circle, fill=red!60!blue!50!white](n6) at (4,-3)  {};
\node [circle, fill=gray!50] (n7) at (0,0)  {};
\node [circle, fill=gray!50] (n8) at (-1.5,-1.8)  {};
\node [circle, fill=gray!50] (n9) at (-3,.5)  {};
\node [circle, fill=gray!50] (n10) at (1.5,-1.8)  {};
\node [circle, fill=gray!50] (n11) at (3,.5)  {};
\draw (n0) -- (n7) [green!40!black]  node [midway, right, black] {};
\draw (n7) -- (n10) [green!40!black]  node [midway,left, black  ] {};
\draw (n7) -- (n8) [green!40!black] node [midway, right, black ] {};
\draw (n8) -- (n9) [green!40!black] node [midway,left , black ] {};
\draw (n8) -- (n3) [green!40!black] node [near end,right , black ] {};
\draw (n9) -- (n2) [green!40!black] node [near end,right , black ] {};
\draw (n9) -- (n1) [green!40!black] node [near end,right , black ] {};
\draw (n10) -- (n4) [green!40!black]  node [midway, black] {};
\draw (n10) -- (n11) [green!40!black]  node [midway, black] {};
\draw (n5) -- (n11) [green!40!black]  node [midway, black] {};
\draw (n6) -- (n11) [green!40!black]  node [midway, black] {};
\end{tikzpicture}} \end{tabular}
& $42$ & $1$ &  $61$ \\
 \hline

 \begin{tabular}{c}\scalebox{0.4}{\begin{tikzpicture}
  [scale=.6,auto=left,very thick,-]
  \node [] (nz) at (3,5.6)  {};
\node [circle, fill=red!60!blue!50!white] (n0) at (-5,5)  {};
\node [circle, fill=red!60!blue!50!white] (n1) at (-8,-1)  {};
\node [circle, fill=red!60!blue!50!white] (n2) at (-7,-5) {};
\node [circle, fill=red!60!blue!50!white](n3) at (3,-5)  {};
\node [circle, fill=red!60!blue!50!white] (n4) at (4,-1) {};
\node [circle, fill=red!60!blue!50!white](n5) at (-2,5)  {};
\node [circle, fill=red!60!blue!50!white](n6) at (1,5)  {};
\node [circle, fill=gray!50] (n7) at (-2,2)  {};
\node [circle, fill=gray!50] (n8) at (-2,0)  {};
\node [circle, fill=gray!50] (n9) at (-5,-2)  {};
\node [circle, fill=gray!50] (n10) at (1,-2)  {};
\draw (n0) -- (n7) [green!40!black]  node [midway, right, black] {};
\draw (n7) -- (n6) [green!40!black]  node [midway,left, black  ] {};
\draw (n7) -- (n5) [green!40!black] node [midway, right, black ] {};
\draw (n7) -- (n8) [green!40!black] node [midway,left , black ] {};
\draw (n8) -- (n9) [green!40!black] node [near end,right , black ] {};
\draw (n8) -- (n10) [green!40!black] node [near end,right , black ] {};
\draw (n10) -- (n3) [green!40!black] node [near end,right , black ] {};
\draw (n10) -- (n4) [green!40!black]  node [midway, black] {};
\draw (n9) -- (n1) [green!40!black]  node [midway, black] {};
\draw (n9) -- (n2) [green!40!black]  node [midway, black] {};
\end{tikzpicture}} \end{tabular}
& $61$ & $4$ &  $99$ &
  \begin{tabular}{c}\scalebox{0.4}{\begin{tikzpicture}
  [scale=.6,auto=left,very thick,-]
\node [circle, fill=red!60!blue!50!white] (n0) at (0,4)  {};
\node [circle, fill=red!60!blue!50!white] (n1) at (-4,4)  {};
\node [circle, fill=red!60!blue!50!white] (n2) at (-5,1) {};
\node [circle, fill=red!60!blue!50!white](n3) at (-5,-1)  {};
\node [circle, fill=red!60!blue!50!white] (n4) at (-4,-4) {};
\node [circle, fill=red!60!blue!50!white](n5) at (4,3.5)  {};
\node [circle, fill=red!60!blue!50!white](n6) at (4,-3.5)  {};
\node [circle, fill=gray!50] (n7) at (0,0)  {};
\node [circle, fill=gray!50] (n8) at (-1.5,0)  {};
\node [circle, fill=gray!50] (n9) at (-3,1.5)  {};
\node [circle, fill=gray!50] (n10) at (-3,-1.5)  {};
\node [circle, fill=gray!50] (n11) at (1.5,0)  {};
\draw (n0) -- (n7) [green!40!black]  node [midway, right, black] {};
\draw (n7) -- (n11) [green!40!black]  node [midway,left, black  ] {};
\draw (n7) -- (n8) [green!40!black] node [midway, right, black ] {};
\draw (n8) -- (n9) [green!40!black] node [midway,left , black ] {};
\draw (n8) -- (n10) [green!40!black] node [near end,right , black ] {};
\draw (n9) -- (n2) [green!40!black] node [near end,right , black ] {};
\draw (n9) -- (n1) [green!40!black] node [near end,right , black ] {};
\draw (n10) -- (n4) [green!40!black]  node [midway, black] {};
\draw (n10) -- (n3) [green!40!black]  node [midway, black] {};
\draw (n5) -- (n11) [green!40!black]  node [midway, black] {};
\draw (n6) -- (n11) [green!40!black]  node [midway, black] {};
\end{tikzpicture}} \end{tabular}
& $53$ & $1$ &  $61$ \\
 \hline
 
     \end{tabular}}
        \caption{Data of BMT models on phylogenetic trees with  $7$ leaves. 
        }
        \label{fig:equivalence_root}
\end{figure}

Computational results show that our formula does not generalize to other BMT models. 
For example, the table in \Cref{fig:equivalence_root} lists all possible non-star tree topologies in $7$ leaves, the degree of the vanishing ideal $\mathcal{I}_T$ (deg), reciprocal maximum likelihood degree (rmld), and the maximum likelihood degree (mld). 
{Since Theorem \ref{thm:main-A} implies that any choice of root results in the same the ML-degree, we do not specify a root for the trees in \Cref{fig:equivalence_root}.}
We computed the  ML-degrees using the software \texttt{HomotopyContinuation.jl} and \texttt{LinearCovarianceModels.jl} \cite{HomotopyContinuation, sturmfels2020estimating}. {Since the number of indeterminates in the likelihood equations grows linearly in the number of leaves, it is also possible to use these strategies to compute the ML-degree for some larger trees fairly quickly. However, as the score equations become more complicated and the number of solutions grows, it becomes more likely that a standard application of these methods will not find all solutions since the paths used for homotopy continuation may cross or be truncated.}

Unlike in the case of the reciprocal ML-degree (see \cite{boege2021reciprocal}), there is not an obvious way to extend the formula for the ML-degree of a star tree to a formula for trees of any topology. The methods that we use in Section \ref{sec:StarTree}  also do not directly extend to arbitrary trees. {When $\T$ is not a star tree,} the common vanishing locus of $\det(K_\T(\theta))$ and the score equations with denominators cleared may be positive dimensional, so one cannot directly apply B\'ezout's theorem. However, we are hopeful that our formulas for $\det(K_\T(\theta))$ and the degree of the path map will be useful in future approaches to this problem.

\section*{Acknowledgements}
\noindent We thank David Speyer for the many insights he shared throughout the course of this project. We also thank Taylor Brysiewicz, Piotr Zwiernik, 
Maximilian Wiesmann, and Leonie Kayser for their helpful comments. { We are grateful to the anonymous reviewers for their careful reading of our manuscript and their detailed feedback.
Thanks especially to Leonie Kayser and an anonymous reviewer for their suggestions to clarify the arguments in Lemmas 5.6 and 5.7.}

This work started at the ``Algebra of phylogenetic networks" workshop at the University of Hawai'i at Mānoa from May 23 - 27, 2022, supported by the National Science Foundation (NSF) under grant DMS-1945584.
JIC and AM were supported by the St John's College Visiting Researcher programme.
JIC was supported by the L'Or\'eal-UNESCO For Women in Science UK and Ireland Rising Talent Award in Mathematics and Computer Science.
SC was supported by the Graduate Research Fellowship under Grant No. DGE-1841052 and by the NSF under Grant No. 1855135. AM was partially supported by the US NSF under Grant No. DMS-2306672 and a Simons Early Career Travel Grant. 
IN was supported by the NSF Grant No. DMS-1945584.
Part of this research was performed while the authors were visiting the Institute for Mathematical and Statistical Innovation (IMSI), which is supported by the NSF Grant No. DMS-1929348.

\bibliographystyle{plain}
\bibliography{references}

\end{document}